\documentclass[reqno]{amsart}
\usepackage{amsfonts}%
\usepackage{mathrsfs}
\usepackage{cite}
\usepackage{color}
\allowdisplaybreaks
\date{\today}
 
\theoremstyle{plain}
\newtheorem{thm}{Theorem}[section]

\newtheorem{prop}[thm]{Proposition}
\theoremstyle{definition}

\theoremstyle{remark}
\newtheorem{rem}{Remark}[section]
\numberwithin{equation}{section}

\renewcommand{\u}{{\bf u}}
\renewcommand{\v}{{\bf v}}

\renewcommand{\H}{{\bf H}}

\newcommand{\dv}{{\rm div }}
\newcommand{\cu}{{\rm curl\, }}

\begin{document}

\title[Incompressible limit of MHD equations]
{Incompressible  limit of the compressible
magnetohydrodynamic equations with periodic boundary conditions}

\author{Song Jiang}
\address{LCP, Institute of Applied Physics and Computational Mathematics, P.O.
 Box 8009, Beijing 100088, P.R. China}
 \email{jiang@iapcm.ac.cn}

\author{Qiangchang Ju}
\address{Institute of Applied Physics and Computational Mathematics, P.O.
 Box 8009-28, Beijing 100088, P.R. China}
 \email{qiangchang\_ju@yahoo.com}

 \author [Fucai Li]{Fucai Li$^*$}%
\thanks{$^*$Corresponding  author}
\address{Department of Mathematics, Nanjing University, Nanjing
 210093, P.R. China}
 \email{fli@nju.edu.cn}



\begin{abstract}
This paper is concerned with the incompressible limit of the
compressible magnetohydrodynamic equations with periodic boundary conditions.
It is rigorously shown that the weak solutions of
the compressible magnetohydrodynamic  equations converge to the strong solution of
 the viscous or inviscid incompressible magnetohydrodynamic
 equations as long as the latter exists both for the well-prepared initial data and general
 initial data. Furthermore, the convergence rates are also obtained
 in the case of the well-prepared initial data.
\end{abstract}

 \maketitle

\section{Introduction}

Magnetohydrodynamics (MHD) studies  the  dynamics of compressible
quasineutrally ionized fluids under the influence of electromagnetic
fields. The applications of magnetohydrodynamics cover a very wide
range of physical objects, from liquid metals to cosmic plasmas.
The
compressible viscous MHD equations in the isentropic case take the
form (see, e.g., \cite{KL,LL,PD})
\begin{align}
&\partial_t \tilde{\rho}  +\dv(\tilde{\rho}\tilde{\u})=0, \label{a1a} \\
&\partial_t
(\tilde{\rho}\tilde{\u})+\dv(\tilde{\rho}\tilde{\u}\otimes\tilde{\u})+\nabla
\tilde{P}  =(\cu\tilde{ \H})\times \tilde{\H}+\tilde{\mu}\Delta\tilde{\u}
+(\tilde{\mu}+\tilde{\lambda})\nabla(\dv\tilde{\u}), \label{a1b} \\
&\partial_t \tilde{\H}
-\cu(\tilde{\u}\times\tilde{\H})=-\cu(\tilde{\nu}\,
\cu\tilde{\H}),\quad \dv\tilde{\H}=0.\label{a1c}
\end{align}
Here $x\in \mathbb{T}^d$, a torus in $\mathbb{R}^d$, $d=2 $ or $3$,  $ t>0$,  the unknowns
$\tilde{\rho}$ denotes the density,
$\tilde{\u}=(\tilde{u}_1,\dots,\tilde{u}_d)\in
 \mathbb{R}^d$ the velocity, and $\tilde{\H}=(\tilde{H}_1,\dots, \tilde{H}_d)\in  \mathbb{R}^d$ the
magnetic field, respectively. The constants  $\tilde{\mu}$ and
$\tilde{\lambda}$ are the shear and bulk viscosity coefficients of
the flow, respectively, satisfying $\tilde{\mu}>0$ and
 $2\tilde{\mu}+d\tilde{\lambda}>0$; the constant $\tilde{\nu}>0$ is the magnetic
diffusivity acting as a magnetic diffusion coefficient of the
magnetic field. $\tilde{P}(\tilde{\rho})$ is the pressure-density
function and here we consider the case
\begin{equation}
\tilde{P}(\tilde{\rho})=a \tilde{\rho}^\gamma,
\label{aad}
\end{equation}
where  $a>0 $ and $ \gamma >1$ are constants.

The well-posedness of the Cauchy problem and initial boundary value problems
for (\ref{a1a})-(\ref{a1c}) has been investigated recently.
The  global existence of weak solutions to the compressible
MHD equations with general initial data was obtained by Hu and Wang~\cite{HW1,HW2}
(also see \cite{FY2} on ``variational solutions'').
  From the physical point of view, one can formally derive the
incompressible models from the compressible ones when the
Mach number goes to zero and the density becomes almost constant.
Based on this observation, Hu and Wang \cite{HW3} proved the convergence of the
weak solutions of the compressible MHD equations \eqref{a1a}-\eqref{a1c}
 to a weak solution of the viscous incompressible  MHD equations.
Jiang, Ju and Li \cite{JJL} obtained the convergence towards the
strong solution of the ideal incompressible MHD equations in the whole space
by using the dispersion property of the wave equation
if both the shear viscosity and the magnetic diffusion coefficients go to zero.

In this paper, we shall extend the results on the Cauchy problem in \cite{JJL} to
the periodic case. First, we consider the well-prepared initial
data for which the oscillations will never appear. We will
rigorously show  the weak solutions of the
compressible MHD equations converge to the strong
solution of the ideal incompressible MHD equations in the periodic
domain if both the shear viscosity and the magnetic
diffusion coefficients go to zero, as well as
to the strong solution of the viscous incompressible  MHD
equations. Furthermore, we shall also give the rates of convergence
  which are not obtained in \cite{HW3, JJL}. Secondly, we
consider the case of general initial data. For this case the oscillations (acoustic waves)
will appear. Comparing with \cite{JJL} where the Cauchy problem was dealt with,
the acoustic waves in the current situation will lose the dispersion
property and will interact each other. Thus, here we have to
impose more regular conditions than $L^2$ on the initial data to control
the oscillating parts. In addition, we have to assume that the Sobolev norm of the
oscillating parts is comparable to the magnetic diffusion coefficient in order to
deal with the general initial data. We will rigorously prove the convergence of the weak
solutions of the compressible MHD equations to the strong solution of the incompressible
MHD equations, as well as to the strong solution of the partial viscous
incompressible MHD equations.

To begin our argument, we first give some formal analysis.
Formally, by utilizing the identity
$$
\nabla(|\tilde{\H}|^2)=2(\tilde{\H}\cdot \nabla
)\tilde{\H}+2\tilde{\H}\times \cu \tilde{\H},
$$
we can rewrite the momentum equation \eqref{a1b}
  as
\begin{equation}\label{a2a}
\partial_t (\tilde{\rho}\tilde{\u})+\dv(\tilde{\rho}\tilde{\u}\otimes\tilde{\u})
+\nabla \tilde{P} =(\tilde{\H}\cdot \nabla
  )\tilde{\H}-\frac12\nabla(|\tilde{\H}|^2)+\tilde{\mu}\Delta\tilde{\u}
  +(\tilde{\mu}+\tilde{\lambda})\nabla(\dv\tilde{\u}).
\end{equation}
By the identities
$$
\cu\cu \tilde{\H}=  \nabla\,\dv \tilde{\H}-\Delta \tilde{\H}$$ and
\begin{equation*}
 \cu( \tilde{ \u}\times\tilde{\H}) =
 \tilde{\u} (\dv \tilde{\H})  -  \tilde{\H}  (\dv\tilde{\u})
 + (\tilde{\H}\cdot \nabla)\tilde{\u} - (\tilde{\u}\cdot \nabla)\tilde{\H},
\end{equation*}
together with the constraint $\dv\tilde{\H}=0$,  the  magnetic field equation
\eqref{a1c} can be   expressed as
\begin{equation}\label{a2b}
 \partial_t\tilde{\H}   +
 (\dv \tilde{\u})\tilde{ \H}+ (\tilde{\u} \cdot \nabla)\tilde{\H}
 - (\tilde{\H}\cdot \nabla)\tilde{\u}= \tilde{\nu }  \Delta \tilde{\H}.
\end{equation}
 We  introduce the  scaling
\begin{equation*}
  \tilde{\rho} (x,  t  )=\rho^\epsilon (x, \epsilon t), \quad
 \tilde{\u} (x, t )=\epsilon \u^\epsilon(x,\epsilon t), \quad
  \tilde{\H} (x, t )=\epsilon \H^\epsilon(x,\epsilon t)
 \end{equation*}
and  assume that the viscosity coefficients $\tilde{\mu}$,
$\tilde{\xi}$, and $\tilde{\nu}$
  are   small constants and scaled like
\begin{equation}\label{pa}
  \tilde{\mu}=\epsilon \mu^\epsilon, \quad \tilde{\lambda}=\epsilon \lambda^\epsilon,
  \quad \tilde{\nu}=\epsilon \nu^\epsilon,
\end{equation}
where $\epsilon\in (0,1)$ is a small parameter and the normalized
coefficients $\mu^\epsilon$, $\lambda^\epsilon$,   and
$\nu^\epsilon$ satisfy  $\mu^\epsilon>0$,
$2\mu^\epsilon+d\lambda^\epsilon>0$, and $\nu^\epsilon>0$.

 With the preceding scalings and the pressure  function \eqref{aad}, the compressible
MHD equations \eqref{a1a}, \eqref{a2a}, and \eqref{a2b} take the
form
\begin{align}
 & \partial_t {\rho}^\epsilon
 +\text{div}({\rho}^\epsilon{\u}^\epsilon)=0,  \label{a2i}\\
&  \partial_t
 ({\rho}^\epsilon{\u}^\epsilon)+\text{div}({\rho}^\epsilon{\u}^\epsilon\otimes
 {\u}^\epsilon)+\frac{a \nabla (\rho^\epsilon)^\gamma
 }{\epsilon^2}
  = ({\H}^\epsilon\cdot \nabla
  ){\H}^\epsilon-\frac12\nabla(|{\H}^\epsilon|^2)\nonumber\\
  &\qquad \qquad\qquad\qquad\qquad\qquad\qquad\qquad\quad
   +{\mu}^\epsilon\Delta{\u}^\epsilon+({\mu}^\epsilon+{\lambda}^\epsilon)
   \nabla(\dv {\u}^\epsilon),   \label{a2j}\\
 &  \partial_t {\H}^\epsilon   +
 (\dv  {\u}^\epsilon) { \H}^\epsilon+ ( {\u}^\epsilon \cdot \nabla) {\H}^\epsilon
 - ( {\H}^\epsilon\cdot \nabla) {\u}^\epsilon=  {\nu }^\epsilon
 \Delta  {\H}^\epsilon,\ \ \dv {\H}^\epsilon=0.\label{a2k}
 \end{align}
Moreover, by replacing $\epsilon$ by $\sqrt{a\gamma}\;\epsilon$, we can always
assume $a=1/\gamma$.

Now, we investigate the incompressible limit of the compressible MHD
equations \eqref{a2i}-\eqref{a2k}. Formally let   $\epsilon
\rightarrow 0$ in the equations \eqref{a2i}-\eqref{a2k}, then  we
obtain from the momentum equation \eqref{a2j} that $\rho^\epsilon $
converges to some function  $\bar \rho(t)\geq 0$. If we further
assume that the initial datum $\rho^\epsilon_0$
 is of order $1+O(\epsilon)$ (this can be guaranteed by the
initial energy bound \eqref{bd} below), then we can expect that
$\bar\rho=1$. Thus, the continuity equation \eqref{a2i} gives ${\rm
div}\,{\u}=0$. Furthermore, using the assumption
\begin{equation}\label{pam}
  \mu^\epsilon\rightarrow 0,    \quad   \nu^\epsilon \rightarrow 0 \quad \text{as}\quad
\epsilon \rightarrow 0,
\end{equation}
  we obtain the following ideal incompressible MHD equations
(suppose that the limits  ${\u}^\epsilon\rightarrow \u$ and
${\H}^\epsilon\rightarrow \H$ exist)
 \begin{align}
& \partial_t \u+(\u\cdot \nabla)\u -({\H} \cdot \nabla
  ){\H}+\nabla p +\frac12\nabla(|{\H} |^2) =0,   \label{a2l}\\
& \partial_t {\H}    + ( {\u}  \cdot \nabla) {\H}
 - ( {\H} \cdot \nabla) {\u} =  0,  \label{a2m}\\
  &\dv \u=0, \quad   \dv \H=0.\label{a2n}
 \end{align}
 In Section 3 we shall rigorously prove that the weak
solutions of the compressible MHD equations \eqref{a2i}-\eqref{a2k}
 converge to, as  $\epsilon\to 0$,
the strong solution of the ideal incompressible MHD equations \eqref{a2l}-\eqref{a2n}
for the  well-prepared initial data in the time interval where
the strong solution of \eqref{a2l}-\eqref{a2n}
exists. Furthermore, the convergence rates are obtained. To show these results,
since the viscosity coefficients
 go to zero, we lose the spatial compactness property of the velocity and the
magnetic field, and the arguments in \cite{HW3} do not work here.
To overcome such difficulty, we shall carefully exploit the energy arguments.

 Next, if we assume that the shear and the bulk viscosity coefficients and
  the magnetic diffusivity coefficient satisfy
\begin{equation}\label{pamb}
  \mu^\epsilon\rightarrow \mu>0, \quad  \lambda^\epsilon \rightarrow \lambda,
\quad\nu^\epsilon \rightarrow \nu>0 \quad \text{as}\quad \epsilon
\rightarrow 0,
\end{equation}
 then the compressible MHD equations \eqref{a2l}-\eqref{a2n}
formally converges to the incompressible MHD equations (suppose that  the limits  ${\u}^\epsilon\rightarrow \u$ and
${\H}^\epsilon\rightarrow \H$ exist)
 \begin{align}
& \partial_t \u+(\u\cdot \nabla)\u -\mu\Delta \u+\nabla p-({\H} \cdot \nabla
  ){\H} +\frac12\nabla(|{\H} |^2) =0,    \label{a2ll}\\
& \partial_t {\H}  + ( {\u}  \cdot \nabla) {\H}
 - ( {\H} \cdot \nabla) {\u} - \nu\Delta \H=  0, \label{a2mm}\\
  &\dv\,\u=0, \quad   \dv \H=0.\label{a2nn}
 \end{align}
In Sections 3 and 4 we shall prove the convergence to the strong
solution of the incompressible viscous MHD equations \eqref{a2ll}-\eqref{a2nn}
for both the well-prepared and the general initial data.
Furthermore, the convergence rates are also obtained for the
well-prepared initial data. For the general initial data, we shall also show
the convergence to the strong solution of the
partial viscous incompressible MHD equations (that is, $\mu=0$ and
$\nu>0$ in \eqref{a2ll}-\eqref{a2nn}).

There are a lot of studies on the compressible MHD
equations in the literature. Besides the aforementioned results,
the interested reader can see \cite{K,LY} on the global smooth solutions with
small initial data and see \cite{VK,FY1} on the local strong solution with
general initial data. We also mention the work \cite{ZJX} where a MHD model
describing the screw pinch problem in plasma physics was discussed
and the global existence of weak solutions with symmetry was obtained.

\smallskip
  Before ending the introduction, we give the notation used throughout the current
  paper. We denote the space $L^q_2(\mathbb{T}^d)$ by
$$
L^q_2(\mathbb{T}^d)= \{f\in L_{loc}(\mathbb{T}^d):  f  1_{\{|f|\geq
1/2\}}\in L^q,  f  1_{\{|f|\leq 1/2\}}\in L^2\} .
$$
We use the letters $C$ and $C_T$ to  denote various positive constants
independent of $\epsilon$, but $C_T$ may depend on $T$. For
convenience, we denote by $H^r  \equiv H^r (\mathbb{T}^d)$ ($r\in\mathbb{R}$) the
standard Sobolev space.   For any vector field $\mathbf{v}$, we denote by
$P\mathbf{v}$ and $Q\mathbf{v}$ the divergence-free
part and the gradient part of $\mathbf{v}$, respectively. Namely,
$Q\mathbf{v}=\nabla \Delta^{-1}(\text{div}\mathbf{v})$ and
$P\mathbf{v}=\mathbf{v}-Q\mathbf{v}$.

  We state our main results in Section 2 and present the proofs for the well-prepared case
  in Section 3 and the ill-prepared case in Section 4, respectively.

\section{main results}

We first recall the local existence of strong solutions to the ideal
incompressible MHD equations \eqref{a2l}-\eqref{a2n} in the torus $\mathbb{T}^d$.
The proof can be found in \cite{DL72,ST}.

\begin{prop}[\!\cite{DL72,ST}]\label{imhd}
Assume that the initial data $(\u ,\H )|_{t=0}=(\u_0,\H_0)$ satisfy
$\u_0,\H_0\in {H}^s$ ($s>d/2+1$), and $\dv\,\u_0=0$, $\dv\H_0=0$. Then, there exist a
 $T^*\in (0,\infty)$ and a unique solution
 $(\u,\H)\in L^{\infty}([0,T^*),{H}^s)$ to the ideal
incompressible MHD equations \eqref{a2l}-\eqref{a2n}
satisfying, for any $0<T<T^*$, $\dv\,\u =0$, $\dv\H =0$, and
\begin{equation}\label{ba}
\sup_{0\le t\le T}\!\big\{||(\u,\H)(t)||_{H^s}
+||(\partial_t\u,\partial_t\H)(t)||_{H^{s-1}}+ ||\nabla p(t)||_{H^{s-1}}\big\} \le C_T.
\end{equation}
\end{prop}
\begin{rem}
The local existence of strong solutions to the incompressible viscous MHD
equations \eqref{a2ll}-\eqref{a2nn} was also established  in \cite{DL72,ST}.
\end{rem}

We prescribe the initial conditions to the compressible
MHD equations   \eqref{a2i}-\eqref{a2k} as
\begin{equation}\label{bb}
  \rho^\epsilon|_{t=0}=\rho^\epsilon_0(x), \quad \rho^\epsilon
 \u^\epsilon|_{t=0}=\rho^\epsilon_0(x)\u^\epsilon_0(x)\equiv \mathbf{m}^\epsilon_0(x),
 \quad \H^\epsilon|_{t=0}=\H^\epsilon_0(x),
   \end{equation}
and assume that
\begin{equation}\label{bc}
 \rho_0^\epsilon \geq 0,\,  \rho^\epsilon_0\in L^\gamma,\,
 \rho^\epsilon_0|\u^\epsilon_0|^2\in L^1,\,\H^\epsilon_0\in L^2,\, \dv\H^\epsilon_0=0,\,
  \mathbf{m}^\epsilon_0=0\;\;\text{ for a.e. }\;\; {\rho^\epsilon_0=0}.
\end{equation}
Moreover,  we assume that the initial data also satisfy the following uniform bound
  \begin{equation}\label{bd}
     \int_{\mathbb{T}^d}
 \Big[\frac12\rho^\epsilon_0|\u^\epsilon_0|^2+\frac12|\H^\epsilon_0|^2
 +\frac{a}{\epsilon^2{(\gamma-1)}}\big((\rho^\epsilon_0)^\gamma -1-\gamma
 (\rho^\epsilon_0-1)\big)
    \Big] dx \leq C.
  \end{equation}
The initial energy inequality \eqref{bd} implies that
$\rho^\epsilon_0$ is of order $1+O(\epsilon)$.

Under the above assumptions, it was proved in \cite{HW1} that the compressible
MHD equations \eqref{a2i}-\eqref{a2k} with initial data \eqref{bb}-\eqref{bd}
has a global weak solution. More precisely, we have

\begin{prop}[\cite{HW1}]\label{cmhd}
 Let $\gamma>d/2$. Suppose that the initial data
 $(\rho^\epsilon_0,\u^\epsilon_0,\H^\epsilon_0)$ satisfy the assumptions \eqref{bc} and \eqref{bd}.
 Then the compressible MHD equations
 \eqref{a2i}-\eqref{a2k} with  the initial data \eqref{bb} enjoy
  at least one global weak solution $(\rho^\epsilon, \u^\epsilon, \H^\epsilon)$
 satisfying
\begin{enumerate}
  \item  $ \rho^\epsilon\in L^\infty(0,\infty;L^\gamma)\cap
 C([0,\infty),L^r)$ for all $1\leq r< \gamma$, $\rho^\epsilon
 |\u^\epsilon|^2\in
  L^\infty(0,\infty; L^1)$, $\H^\epsilon
\in   L^\infty(0,\infty; L^2)$,
 and  $\u^\epsilon \in L^2(0,T; H^1)$,
  $\rho^\epsilon\u^\epsilon\in C([0,T],   L^{\frac{2\gamma}{\gamma+1}}_{weak})$,
   $\H^\epsilon \in L^2(0,T; H^1)\cap  C([0,T],L^{\frac{2\gamma}{\gamma+1}}_{weak})$
   for all $T\in (0,\infty)$;

  \item  the energy inequality
  \begin{equation}\label{be}
   \mathcal{E}^\epsilon(t)+\int^t_0  \mathcal{D}^\epsilon(s)ds\leq
 \mathcal{E}^\epsilon(0)
  \end{equation}
  holds with the finite total energy
  \begin{equation}\label{bf}
 \qquad\qquad  \mathcal{E}^\epsilon(t)\equiv \int_{\mathbb{T}^d}
 \Big[\frac12\rho^\epsilon|\u^\epsilon|^2+\frac12|\H^\epsilon|^2
        +\frac{a}{\epsilon^2{(\gamma-1)}}\big((\rho^\epsilon)^\gamma
 -1-\gamma (\rho^\epsilon-1)\big)
    \Big](t)
   \end{equation}
   and the dissipation energy
\begin{equation}\label{bff}
\mathcal{D}^\epsilon(t)\equiv \int_{\mathbb{T}^d}\big[ \mu^\epsilon|\nabla
 \u^\epsilon|^2
 +(\mu^\epsilon+\lambda^\epsilon) |{\rm div}\u^\epsilon|^2 + \nu^\epsilon |\nabla
 \H^\epsilon|^2\big](t);
\end{equation}

  \item the continuity equation is satisfied in the sense of
 renormalized solutions, i.e.,
  \begin{equation} \label{cnsbg}
    \partial_t b(\rho^\epsilon)+{ \rm div}
 (b(\rho^\epsilon)\u^\epsilon)+\big(b'(\rho^\epsilon)\rho^\epsilon-b(\rho^\epsilon)\big)
 {\rm div}\u^\epsilon =0
  \end{equation}
for any $b\in C^1(\mathbb{T})$ such that $b'(z)$ is constant for $z$ large enough;
  \item the equations \eqref{a2i}-\eqref{a2k} hold in
 $\mathcal{D}'( \mathbb{T}^d\times(0,\infty))$.
\end{enumerate}
\end{prop}

 The main results of this paper can be stated as follows.
\begin{thm}\label{MRa} Let $s>{d}/{2}+2$ and $\mu^\epsilon+\lambda^\epsilon>0$.
 Suppose that the initial data $(\rho^\epsilon_0,\u^\epsilon_0,\H_0^\epsilon)$
 satisfy the conditions presented in Proposition \ref{cmhd}.
 Assume further that
\begin{align}
&\int_{\mathbb{T}^d}|\rho_0^\epsilon-1|^21_{(|\rho_0^\epsilon-1|\leq
\delta)}\;dx+\int_{\mathbb{T}^d}|\rho_0^\epsilon-1|^\gamma1_{(|\rho_0^\epsilon-1|>
\delta)}\;dx \leq C\epsilon^2,\label{icd1}\\
& \qquad
||\sqrt{\rho^\epsilon_0}\mathbf{u}_0^\epsilon-\mathbf{u}_0||^2_{L^2(\mathbb{T}^d)}
 \leq C\epsilon, \quad ||\mathbf{H}_0^\epsilon-\mathbf{H}_0||^2_{L^2(\mathbb{T}^d)}
\leq C\epsilon \label{icd2}
\end{align}
for any $\delta\in (0,1)$, where $\mathbf{u}_0$ and $\mathbf{H}_0$ are defined
in Proposition \ref{imhd}. We assume that the shear viscosity $\mu^\epsilon$
and the magnetic diffusion coefficient $\nu^\epsilon $ satisfy
\begin{align}\label{abc}
  \mu^\epsilon =\epsilon^\alpha,\quad \nu^\epsilon=\epsilon^\beta
\end{align}
for some constants $\alpha , \beta>0$ satisfying $0<\alpha+\beta <2$.
  Let $(\u,\H)$ be the smooth
solution to the ideal incompressible MHD equations
\eqref{a2l}-\eqref{a2n} defined on $[0,T^*)$ with
$(\u ,\H )|_{t=0}=(\u_0,\H_0)$. Then, for any $0<T<T^*$, the
global weak solution $(\rho^\epsilon, \u^\epsilon, \H^\epsilon)$ of
the compressible MHD equations \eqref{a2i}-\eqref{a2k} established
in Proposition \ref{cmhd} satisfies
\begin{align}
&\int_{\mathbb{T}^d}|\rho^\epsilon-1|^21_{(|\rho^\epsilon-1|\leq
\delta)}\;dx+\int_{\mathbb{T}^d}|\rho^\epsilon-1|^\gamma1_{(|\rho^\epsilon-1|>
\delta)}\;dx \leq C_T\epsilon^2, \label{icd10}\\
& \qquad
||\sqrt{\rho^\epsilon}\mathbf{u}^\epsilon-\mathbf{u}||^2_{L^2(\mathbb{T}^d)}
\leq C_T\epsilon^\sigma, \quad
||\mathbf{H}^\epsilon-\mathbf{H}||^2_{L^2(\mathbb{T}^d)} \leq
C_T\epsilon^\sigma \label{icd20}  \end{align} for any $t\in [0,T]$,
where $\sigma=\min\{\alpha,\beta, 1-(\alpha+\beta)/2 \}.$
\end{thm}

The proof of Theorem \ref{MRa} is based on  the combination  of  the modulated energy method,
motivated by Brenier \cite{B00}, the weak convergence method and the refined energy
analysis. Masmoudi \cite{M01b} made use of
such idea to study the incompressible, inviscid limit
of  the compressible Navier-Stokes equations in both the whole space and the torus.
Comparing with the proof in \cite{M01b}, here we have to overcome the difficulties caused by
the strong coupling of the hydrodynamic motion and  the magnetic field.

Furthermore, we can use an idea similar to that described above to
obtain the convergence of the compressible MHD equations \eqref{a2i}-\eqref{a2k}
to the incompressible viscous MHD equations \eqref{a2ll}-\eqref{a2nn}.
In fact,  we have the following result.

\begin{thm}\label{MRb}  Let $s>{d}/{2}+2$ and $\mu^\epsilon+\lambda^\epsilon>0$.
  Suppose that the initial data $(\rho^\epsilon_0,\u^\epsilon_0,\H_0^\epsilon)$
satisfy the conditions presented   in Proposition \ref{cmhd}. Assume further that
\begin{align}
&\int_{\mathbb{T}^d}|\rho_0^\epsilon-1|^21_{(|\rho_0^\epsilon-1|\leq
\delta)}\;dx+\int_{\mathbb{T}^d}|\rho_0^\epsilon-1|^\gamma1_{(|\rho_0^\epsilon-1|>
\delta)}\;dx \leq C\epsilon^2,\label{icd11}\\
& \qquad ||\sqrt{\rho^\epsilon_0}\mathbf{u}_0^\epsilon
-\mathbf{u}_0||_{L^2(\mathbb{T}^d)}^2 \leq C\epsilon, \quad
||\mathbf{H}_0^\epsilon-\mathbf{H}_0||_{L^2(\mathbb{T}^d)}^2 \leq
C\epsilon \label{icd22}
\end{align}
for any $\delta\in (0,1)$ and for some $\mathbf{u}_0,\mathbf{H}_0\in
H^s(\mathbb{T}^d)$, satisfying $\dv \u_0=0,   \dv \H_0 =0.$ We also
assume that the shear viscosity $\mu^\epsilon $ and the magnetic
diffusion coefficient $\nu^\epsilon $ satisfy \eqref{pamb}.
 Let $(\u,\H)$ be the smooth
solution to the  incompressible MHD equations \eqref{a2ll}-\eqref{a2nn}
 with $(\u ,\H )|_{t=0}=(\u_0,\H_0)$.
 Then, for any $0<T<T^{**}$ ( $ T^{**}$ is the maximal time of existence
for \eqref{a2ll}-\eqref{a2nn}), the
global weak solution $(\rho^\epsilon, \u^\epsilon, \H^\epsilon)$ of the
compressible MHD equations
\eqref{a2i}-\eqref{a2k}
established in Proposition \ref{cmhd} satisfies that
 $\nabla \u^\epsilon$ and  $\nabla \H^\epsilon$ converge strongly to $\nabla \u$
 and $\nabla \H$ in $L^2 (0,T;L^2(\mathbb{T}^d))$, respectively.
 Moreover, for any $t\in [0,T]$, we have
\begin{align}
&\int_{\mathbb{T}^d}|\rho^\epsilon-1|^21_{(|\rho^\epsilon-1|\leq
\delta)}\;dx+\int_{\mathbb{T}^d}|\rho^\epsilon-1|^\gamma1_{(|\rho^\epsilon-1|>
\delta)}\;dx \leq C_T\epsilon^2,\label{icd110}\\
&
||\sqrt{\rho^\epsilon}\mathbf{u}^\epsilon-\mathbf{u}||^2_{L^2(\mathbb{T}^d)}
+||\mathbf{H}^\epsilon-\mathbf{H}||^2_{L^2(\mathbb{T}^d)} \leq
C_T\frac{\epsilon}{\sqrt{\mu\nu}}+C_T\Big(\frac{|\mu^\epsilon-\mu|}{\sqrt{\mu}}
+\frac{|\nu^\epsilon-\nu|}{\sqrt{\nu}}\Big).\label{icd220}
\end{align}
\end{thm}

To show Theorem \ref{MRb}, besides the techniques mentioned above, we have to employ
a new technique, that is, to modulate both the total energy and the partial dissipative
energy simultaneously. Moreover, the dissipative effect of the viscous terms is also
carefully exploited to obtain the desired results.

\begin{rem}
Comparing with Theorem \ref{MRa}, we have gotten the better convergence rates \eqref{icd220}
than (\ref{icd20}) when the
shear viscosity and the magnetic diffusion coefficient tend to some positive constants.
\end{rem}

\medskip
Some results in Theorem \ref{MRb} can be extended to the case of
general initial data. More precisely, we shall obtain
  the convergence of the compressible MHD equations \eqref{a2i}-\eqref{a2k}
 to  the  incompressible MHD
 equations  \eqref{a2ll}-\eqref{a2nn} for the general initial data
 under the conditions that the oscillating parts  of the initial data have
higher regularity and the Sobolev norm of the oscillation parts
 is comparable to the magnetic diffusion coefficient.
This implies that the influence of oscillations
   on the magnetic field can be balanced by the diffusive effect of the magnetic field,
   which is one of the new ingredients in our paper.

   To describe the result, we write $\rho^\epsilon=1+\epsilon
\varphi^\epsilon$ and denote
\begin{equation*}
  \Pi^\epsilon(x,t)=\frac{1}{\epsilon}\sqrt{\frac{2a}{\gamma-1}
  \big((\rho^\epsilon)^\gamma-1-\gamma(\rho^\epsilon-1)\big)}.
\end{equation*}
We will use the above approximation $\Pi^\epsilon(x,t)$ instead of
$\varphi^\epsilon$, since we can not obtain any bound for
$\varphi^\epsilon$ in $L^\infty(0,T;L^2)$ directly if $\gamma <2$.

\begin{thm}\label{MRc}  Let $s>2+{d}/{2}$ and $2\mu+d\lambda>0$.
  Suppose that the initial data $(\rho^\epsilon_0,\u^\epsilon_0,\H_0^\epsilon)$
   satisfy the conditions presented   in Proposition \ref{cmhd}.
  Moreover, we assume that
$\sqrt{\rho^\epsilon_0}\u^\epsilon_0$ converges strongly in $L^2$ to
some $\tilde{\u}_0$ satisfying $Q\tilde{\u}_0\in H^{s-1}$,
$\H^\epsilon_0$ converges strongly in $L^2$ to some $\H_0$ with
$\int_{\mathbb{T}^d} \H_0(x)dx=0$,
$\Pi^\epsilon|_{t=0}=\Pi^\epsilon_0$ converges strongly in $L^2$ to
some $\varphi_0\in H^{s-1}$, and
 \begin{align}\label{Qu}
\|\varphi_0\|_{H^2}+||Q\tilde{ \u}_0||_{H^2}\leq c_0 \nu
\end{align}
 for some   constant $c_0>0$. Let $(\u,\H)$ be the smooth
solution to the  incompressible MHD equations \eqref{a2ll}-\eqref{a2nn}
 with $(\u ,\H )|_{t=0}=(\u_0,\H_0)\in H^s(\mathbb{T}^d)$ satisfying
 $\u_0=P\tilde{ \u}_0$ and $\dv \H_0=0$.
 Then, for any $0<T<T^{**}$ ( $ T^{**}$ is the maximal time of existence
for \eqref{a2ll}-\eqref{a2nn}), the
global weak solution $(\rho^\epsilon, \u^\epsilon, \H^\epsilon)$ of the
compressible MHD equations
\eqref{a2i}-\eqref{a2k}
established in Proposition \ref{cmhd} satisfies
\begin{enumerate}
          \item  $\rho^\epsilon$ converges strongly to $1$ in
 $C([0,T],L^\gamma_2(\mathbb{T}^d))$;
          \item  $\nabla \H^\epsilon$ converges strongly to $\nabla \H$
 in $L^2 (0,T;L^2(\mathbb{T}^d))$;
\item $\H^\epsilon$ converges strongly to
 $\H$ in $L^\infty(0,T; L^2(\mathbb{T}^d))$;
          \item $P(\sqrt{\rho^\epsilon}\u^\epsilon)$ converges  strongly
 to $\u$ in $L^\infty (0,T;L^2(\mathbb{T}^d))$;
          \item $\sqrt{\rho^\epsilon}\u^\epsilon$ converges weakly to
          $\u$ in $H^{-1}(0,T; L^2(\mathbb{T}^d))$.
         \end{enumerate}
\end{thm}

By slightly modifying the proof of Theorem \ref{MRc}, we can obtain
the convergence of compressible MHD equations to the partial viscous
incompressible MHD equations when the shear viscosity goes to zero
and the magnetic diffusion coefficient goes to a positive constant.
The partial viscous incompressible MHD equations correspond to the
case of turbulent flow with very high Reynolds  number (where the
viscosity of flow can be ignored, see \cite{LZ}).
\begin{thm}\label{MRd}  Let $s>2+{d}/{2}$.
  Suppose that the conditions  in Theorem \ref{MRc} hold.    Moreover,
we assume that
\begin{align*}
  \nu^\epsilon \rightarrow \nu>0, \quad 2\mu^\epsilon+\lambda^\epsilon
  \rightarrow 2\theta> 0 \quad \text{as} \quad \epsilon \rightarrow 0,
\end{align*}
and $\mu^\epsilon =\epsilon^\alpha $ for some constant $0< \alpha<1$.   Let $(\u,\H)$ be the smooth
solution to the following partially viscous  incompressible MHD equations
 \begin{align*}
& \partial_t \u+(\u\cdot \nabla)\u  +\nabla p-({\H} \cdot \nabla
  ){\H} +\frac12\nabla(|{\H} |^2) =0,     \\
& \partial_t {\H}  + ( {\u}  \cdot \nabla) {\H}
 - ( {\H} \cdot \nabla) {\u} - \nu\Delta \H=  0,  \\
  &\dv \u=0, \quad   \dv \H=0,
 \end{align*}
 with $(\u ,\H )|_{t=0}=(\u_0,\H_0)\in H^s(\mathbb{T}^3)$ satisfying
 $\u_0=P\tilde{ \u}_0$ and $\dv \H_0=0$.
 Then, for any $0<T<T^{**}$ ($ T^{**}$ is the maximal time of existence for
\eqref{a2ll}-\eqref{a2nn}), the
global weak solution $(\rho^\epsilon, \u^\epsilon, \H^\epsilon)$ of
the compressible MHD equations \eqref{a2i}-\eqref{a2k} established
in Proposition \ref{cmhd} satisfies
\begin{enumerate}
          \item  $\rho^\epsilon$ converges strongly to $1$ in
 $C([0,T],L^\gamma_2(\mathbb{T}^d))$;
          \item  $\nabla \H^\epsilon$ converges strongly to $\nabla \H$
 in $L^2 (0,T;L^2(\mathbb{T}^d))$;
\item $\H^\epsilon$ converges strongly to
 $\H$ in $L^\infty(0,T; L^2(\mathbb{T}^d))$;
          \item $P(\sqrt{\rho^\epsilon}\u^\epsilon)$ converges  strongly
 to $\u$ in $L^\infty (0,T;L^2(\mathbb{T}^d))$;
           \item $\sqrt{\rho^\epsilon}\u^\epsilon$ converges weakly to
           $\u$ in $H^{-1}(0,T; L^2(\mathbb{T}^d))$.
         \end{enumerate}
\end{thm}

\begin{rem}
The assumption that $\Pi^\epsilon_0$ converges strongly in $L^2$ to
some $\varphi_0$ in fact implies that $\varphi^\epsilon_0$ converges
strongly  to $\varphi_0$ in $L^\gamma_2$.
\end{rem}

\begin{rem}
When  taking $\H^\epsilon\equiv 0$ in \eqref{a1a}-\eqref{a1c}, the
MHD equations reduce to the classical compressible Navier-Stokes equations.
The low Mach number limit problem of the compressible Navier-Stokes equations
has been investigated extensively, for instance, see \cite{G,H98,LM98,DG99,D02}. The interested
reader  can refer to the survey article \cite{M06} for more related results.
\end{rem}

\begin{rem}
  We point out that our arguments in the present paper
can be applied to the case of $\H^\epsilon\equiv 0$. In this case, we obtain the
convergence of the compressible Navier-Stokes equtions to the incompressible
 Euler or Navier-Stokes equations with general initial data, extending thus the results
 in \cite{LM98,M01b}.
\end{rem}

\section{Proof of  Theorems \ref{MRa} and \ref{MRb}}

In this section, we shall prove our convergence results for the case
of well-prepared initial data  by combining the modulated energy
method, the weak convergence method, and the refined energy
analysis.

\begin{proof}[Proof of Theorem \ref{MRa}]
  We divide the proof into several steps.

\medskip
\emph{Step 1: Basic energy estimates and compact arguments.}

By the assumptions on the initial data  we obtain, from the energy
inequality \eqref{be}, that the total energy
$\mathcal{E}^\epsilon(t)$ has a uniform upper bound  for a.e. $t\in
[0,T]$, $T>0$. This uniform bound
implies that $\rho^\epsilon |\u^\epsilon|^2$ and
$ \big((\rho^\epsilon)^\gamma-1-\gamma(\rho^\epsilon-1)\big)/{\epsilon^2}$
are bounded in $L^\infty(0,T;L^1)$ and $ \H^\epsilon $ is
bounded in $L^\infty(0,T;L^2)$. Using the analysis in \cite{LM98},  we obtain
\begin{equation}\label{L2g}
\int_{\mathbb{T}^d}\frac{1}{\epsilon^2}|\rho^\epsilon-1|^2
1_{\{|\rho^\epsilon-1|\leq \frac12\}}
  +\int_{\mathbb{T}^d}\frac{1}{\epsilon^2}|\rho^\epsilon-1|^\gamma
  1_{\{|\rho^\epsilon-1|\geq \frac12\}}   \leq C,
\end{equation}
which implies \eqref{icd10}
 and
\begin{equation}\label{re}
    \rho^\epsilon \rightarrow    1 \ \
    \text{strongly in}\ \
C([0,T],L^\gamma_2(\mathbb{T}^d)).
\end{equation}
From the results in \cite{LM98}, we know that
$\|\u^\epsilon\|_{L^2_tL^2_x}^2\leq
C+C\|\nabla\u^\epsilon\|_{L^2_tL^2_x}^2$. Furthermore, the fact that
$\rho^\epsilon |\u^\epsilon|^2$ and $ |\H^\epsilon|^2$ are bounded
in $L^\infty(0,T;L^1)$ implies the following convergence (up to the
extraction of a subsequence $\epsilon_n$):
\begin{align*}
& \sqrt{\rho^\epsilon}  \u^\epsilon   \  \text{converges weakly-$\ast$}
\,\, \text{to some}\, \mathbf{J} \, \,\text{in}\,
L^\infty(0,T;L^2(\mathbb{T}^d)),\\
&  \H^\epsilon \    \text{converges weakly-$\ast$}
\,\, \text{to some}\, \mathbf{K} \, \,\text{in}\,
L^\infty(0,T;L^2(\mathbb{T}^d)).
\end{align*}

Thus, to finish our proof, we need to
  show that $\mathbf{J}=\u$ and  $\mathbf{K}=\H$ in some sense and the inequalities \eqref{icd20} hold,
where $(\u,\H)$  is the strong solution to the ideal
incompressible MHD equations \eqref{a2l}-\eqref{a2n}.


\medskip
\emph{Step 2: The modulated energy functional and the uniform estimates.}

We first recall the energy inequality of the compressible
MHD equations \eqref{a2i}-\eqref{a2k}, i.e., for almost all
$t$, there holds
 \begin{align}\label{cei}
 & \frac12\int_{\mathbb{T}^d}\Big[\rho^\epsilon(t)|\u^\epsilon|^2(t)+|\H^\epsilon|^2(t)
    +(\Pi^\epsilon(t))^2
      \Big]
+\mu^\epsilon \int^t_0\! \int_{\mathbb{T}^d}|\nabla \u^\epsilon|^2 \nonumber\\
&\ \  +(\mu^\epsilon+\lambda^\epsilon)\int^t_0\! \int_{\mathbb{T}^d}|{\rm div}\u^\epsilon |^2
+\nu^\epsilon\int^t_0\! \int_{\mathbb{T}^d}|\nabla \H^\epsilon|^2\nonumber\\
& \ \ \   \leq  \frac12\int_{\mathbb{T}^d}\Big[\rho_0^\epsilon|\u_0^\epsilon|^2
+|\H_0^\epsilon|^2+(\Pi^\epsilon_0)^2     \Big] .
 \end{align}
The conservation of energy for the ideal incompressible MHD equations
\eqref{a2l}-\eqref{a2n} reads
\begin{equation}\label{cfi}
  \frac12 \int_{\mathbb{T}^d}\big[|\u|^2(t)+ |\H|^2(t)\big]   =\frac12  \int_{\mathbb{T}^d}
  \big[|\u_0|^2+|\H_0|^2\big].
\end{equation}
Using $\u$ to test the  momentum equation \eqref{a2j}, we obtain
\begin{align}\label{cii}
&  \int_{\mathbb{T}^d}(\rho^\epsilon \u^\epsilon \cdot \u)(t)
+\int^t_0\int_{\mathbb{T}^d}\rho^\epsilon \u^\epsilon \cdot \big[(\u\cdot \nabla)
\u-(\H\cdot \nabla)\H+\nabla p+\frac12\nabla(|\H|^2)\big] \nonumber\\
&  -\int^t_0\int_{\mathbb{T}^d}\big[(\rho^\epsilon \u^\epsilon \otimes
\u^\epsilon)\cdot \nabla \u+   (\H^\epsilon\cdot \nabla)\H^\epsilon \cdot \u
 - \mu^\epsilon  \nabla  \u^\epsilon  \cdot \nabla \u\big] = \int_{\mathbb{T}^d}
  \rho^\epsilon_0 \u^\epsilon_0 \cdot \u_0.
\end{align}
 Similarly, using  $\H$ to test the magnetic field equation \eqref{a2k}, one gets
\begin{align}\label{cj1i}
&  \int_{\mathbb{T}^d}(\H^\epsilon \cdot \H)(t)
+\int^t_0\int_{\mathbb{T}^d}\H^\epsilon \cdot\big[(\u\cdot \nabla)
\H-(\H\cdot \nabla)\u\big] + \nu^\epsilon
\int^t_0\int_{\mathbb{T}^d}
\nabla  \H^\epsilon  \cdot \nabla \H\nonumber\\
&\quad +\int^t_0\int_{\mathbb{T}^d}\big[(\dv  {\u}^\epsilon) { \H}^\epsilon
+ ( {\u}^\epsilon \cdot \nabla) {\H}^\epsilon  - ( {\H}^\epsilon\cdot \nabla)
 {\u}^\epsilon\big]\cdot \H  = \int_{\mathbb{T}^d} \H^\epsilon_0 \cdot \H_0.
\end{align}
Summing up \eqref{cei} and \eqref{cfi}, and inserting
\eqref{cii} and \eqref{cj1i} into the resulting inequality,
 we can deduce the following
inequality by a straightforward computation
\begin{align}\label{cki}
&   \frac{1}{2}\int_{\mathbb{T}^d}\Big\{  |\sqrt{\rho^\epsilon} \u^\epsilon
-\u|^2(t)+ |  \H^\epsilon
-\H|^2(t)
  +(\Pi^\epsilon)^2(t)\Big\}\nonumber\\
&\quad
+ {\mu^\epsilon}  \int^t_0\int_{\mathbb{T}^d}|\nabla \u^\epsilon|^2
+ (\mu^\epsilon+\lambda^\epsilon)\int^t_0\! \int_{\mathbb{T}^d}|{\rm div}
\u^\epsilon |^2+ {\nu^\epsilon} \int^t_0\int_{\mathbb{T}^d}|\nabla \H^\epsilon |^2\nonumber\\
& \leq  {\mu^\epsilon}  \int^t_0\int_{\mathbb{T}^d}\nabla \u^\epsilon\cdot\nabla \u
+       {\nu^\epsilon}  \int^t_0\int_{\mathbb{T}^d}\nabla \H^\epsilon\cdot\nabla \H
           -\int^t_0\int_{\mathbb{T}^d}\rho^\epsilon
\u^\epsilon\cdot [(\H\cdot \nabla) \H)]       \nonumber\\
&\quad  -\int^t_0\int_{\mathbb{T}^d} (\H^\epsilon\cdot \nabla)
\H^\epsilon \cdot\u  +\int^t_0\int_{\mathbb{T}^d}\H^\epsilon
\cdot\big[(\u\cdot \nabla) \H-(\H\cdot \nabla)\u\big] \nonumber\\
&\quad
 +\int^t_0\int_{\mathbb{T}^d}\big[(\dv  {\u}^\epsilon) { \H}^\epsilon+
 ( {\u}^\epsilon \cdot \nabla) {\H}^\epsilon
 - ( {\H}^\epsilon\cdot \nabla) {\u}^\epsilon\big]\cdot \H
 +\frac12\int^t_0\int_{\mathbb{T}^d}\rho^\epsilon \u^\epsilon \cdot \nabla(|\H|^2)
 \nonumber\\
&\quad +\int^t_0\int_{\mathbb{T}^d}\rho^\epsilon \u^\epsilon \cdot\big[ (\u\cdot
\nabla) \u +\nabla p \big]-\int^t_0\int_{\mathbb{T}^d}(\rho^\epsilon
\u^\epsilon\otimes \u^\epsilon)\cdot \nabla \u\nonumber\\
&\quad
+\int_{\mathbb{T}^d}
\big[(\sqrt{\rho^\epsilon}-1)\sqrt{\rho^\epsilon}\u^\epsilon\cdot
\u\big](t) -\int_{\mathbb{T}^d}
\big[(\sqrt{\rho^\epsilon}-1)\sqrt{\rho^\epsilon}\u^\epsilon\cdot
\u\big](0)\nonumber\\
&\quad   +\frac{1}{2}\int_{\mathbb{T}^d}\Big\{  |\sqrt{\rho^\epsilon} \u^\epsilon
-\u|^2(0)+|\H^\epsilon
-\H |^2(0)  +(\Pi^\epsilon_0)^2\Big\}.
\end{align}

We first deal with the right-hand side of the inequality
\eqref{cki}. Denoting  $\mathbf{w}^{\epsilon}=\sqrt{\rho^\epsilon}
\u^\epsilon -\u$ and  $\mathbf{Z}^{\epsilon}= \H^\epsilon -\H$,
integrating by parts, and using the fact that ${\rm div}\,
\H^\epsilon =0,$ $ {\rm div}\, \u =0$ and ${\rm div}\,\H =0$, we
find that
\begin{align}\label{ck1i}
&   -\int^t_0\int_{\mathbb{T}^d}\rho^\epsilon
\u^\epsilon\cdot [(\H\cdot \nabla) \H)]
    -\int^t_0\int_{\mathbb{T}^d}
(\H^\epsilon\cdot \nabla) \H^\epsilon \cdot\u \nonumber\\
& +\int^t_0\int_{\mathbb{T}^d}\H^\epsilon \cdot\big[(\u\cdot \nabla)
\H-(\H\cdot
\nabla)\u\big]+\frac{1}{2}\int^t_0\int_{\mathbb{T}^d}\rho^\epsilon
\u^\epsilon\cdot \nabla(|\H|^2)\nonumber\\
&
 +\int^t_0\int_{\mathbb{T}^d}\big[(\dv  {\u}^\epsilon) { \H}^\epsilon
 + ( {\u}^\epsilon \cdot \nabla) {\H}^\epsilon
 - ( {\H}^\epsilon\cdot \nabla) {\u}^\epsilon\big]\cdot \H  \nonumber\\
 = &  -\int^t_0\int_{\mathbb{T}^d}\rho^\epsilon
\u^\epsilon\cdot [(\H\cdot \nabla) \H)] +\int^t_0\int_{\mathbb{T}^d}
(\H^\epsilon\cdot \nabla) \u \cdot \H^\epsilon \nonumber\\
& +\int^t_0\int_{\mathbb{T}^d}(\u\cdot \nabla)\H\cdot \H^\epsilon
-\int^t_0\int_{\mathbb{T}^d}(\H\cdot \nabla)\u  \cdot \H^\epsilon  \nonumber\\
&  - \int^t_0\int_{\mathbb{T}^d}( {\u}^\epsilon \cdot \nabla)\H \cdot {\H}^\epsilon
 + \int^t_0\int_{\mathbb{T}^d}( {\H}^\epsilon\cdot \nabla) \H \cdot {\u}^\epsilon
 +\frac{1}{2}\int^t_0\int_{\mathbb{T}^d}\rho^\epsilon
\u^\epsilon\cdot \nabla(|\H|^2)  \nonumber\\
 = & \int^t_0\int_{\mathbb{T}^d}(1-\rho^\epsilon)
\u^\epsilon\cdot [(\H\cdot \nabla) \H)]
+\int^t_0\int_{\mathbb{T}^d}[(\H^\epsilon-\H)\cdot \nabla]\u  \cdot (\H^\epsilon-\H)\nonumber\\
&+\int^t_0\int_{\mathbb{T}^d}[(\H^\epsilon-\H)\cdot \nabla]\H  \cdot (\u^\epsilon-\u)
- \int^t_0\int_{\mathbb{T}^d}[(\u^\epsilon-\u)\cdot \nabla]\H  \cdot (\H^\epsilon-\H)\nonumber\\
& +\frac12\int^t_0\int_{\mathbb{T}^d}(\rho^\epsilon-1){\u}^\epsilon \nabla(| \H|^2)\nonumber\\
\leq & \int^t_0\int_{\mathbb{T}^d}(1-\rho^\epsilon)
\u^\epsilon\cdot [(\H\cdot \nabla) \H)]
+\int^t_0    ||\mathbf{Z}^\epsilon(s)||^2_{L^2}|| \nabla \u(s)||_{L^\infty}  ds\nonumber\\
&+  \int^t_0   \big[||\mathbf{w}^{\epsilon}(s)||^2_{L^2}
+ ||\mathbf{Z}^\epsilon(s)||^2_{L^2}\big] || \nabla \H(s)||_{L^\infty}  ds
+\int^t_0\!\int_{\mathbb{T}^d}(\mathbf{Z}^\epsilon \cdot \nabla)\H
\cdot [(1-\sqrt{\rho^\epsilon})\u^\epsilon]
\nonumber\\
&- \int^t_0\!\int_{\mathbb{T}^d}\big\{[(1-\sqrt{\rho^\epsilon})\u^\epsilon
]\cdot \nabla\big\}\H  \cdot   \mathbf{Z}^\epsilon
+\frac12\int^t_0\int_{\mathbb{T}^d}(\rho^\epsilon-1){\u}^\epsilon \nabla(|\H|^2)
\end{align}
and
\begin{align}\label{cli}
 & \int^t_0\int_{\mathbb{T}^d}\big[\rho^\epsilon \u^\epsilon \cdot((\u\cdot
\nabla) \u+\nabla p)\big]-\int^t_0\int_{\mathbb{T}^d}(\rho^\epsilon
\u^\epsilon\otimes \u^\epsilon)\cdot \nabla \u\nonumber\\
    =& -\int^t_0\int_{\mathbb{T}^d}(\mathbf{w}^{\epsilon}\otimes
 \mathbf{w}^{\epsilon})\cdot \nabla \u
       +\int^t_0\int_{\mathbb{T}^d}(\rho^\epsilon-\sqrt{\rho^\epsilon})\u^\epsilon\cdot
 ((\u\cdot \nabla) \u)\nonumber\\
   &
    +\int^t_0\int_{\mathbb{T}^d}[(\sqrt{\rho^\epsilon}\u^\epsilon-\u)\cdot \nabla] \u
 \cdot \mathbf{w}^{\epsilon}
      + \int^t_0 \int_{\mathbb{T}^d}\rho^\epsilon \u^\epsilon \cdot  \nabla p\nonumber\\
    &  -\int^t_0\int_{\mathbb{T}^d}(\sqrt{ \rho^\epsilon}\u^\epsilon- \u) \cdot \nabla
 \big(\frac{|\u|^2}{2}\big).
    \end{align}

Substituting \eqref{ck1i} and \eqref{cli} into the inequality \eqref{cki}, we conclude that
\begin{align}\label{cmi}
 &
 ||\mathbf{w}^{\epsilon}(t)||^2_{L^2}+
 ||\mathbf{Z}^{\epsilon}(t)||^2_{L^2}+||\Pi^\epsilon(t)||^2_{L^2}
  + 2{\mu^\epsilon} \int^t_0\int_{\mathbb{T}^d}|\nabla \u^\epsilon|^2 \nonumber\\
& + 2(\mu^\epsilon+\lambda^\epsilon)\int^t_0\! \int_{\mathbb{T}^d}|{\rm
div}\u^\epsilon |^2 + 2{\nu^\epsilon} \int^t_0\int_{\mathbb{T}^d}|\nabla \H^\epsilon|^2
\nonumber\\
\leq & 2C \int^t_0 (||\mathbf{w}^{\epsilon}(\tau)||^2_{L^2}+
 ||\mathbf{Z}^{\epsilon}(\tau)||^2_{L^2})(||\nabla \u(\tau)||_{L^\infty}+||\nabla
\H(\tau)||_{L^\infty})d\tau\nonumber\\
& +||\mathbf{w}^{\epsilon}(0)||^2_{L^2}+||\mathbf{Z}^{\epsilon}(0)||^2_{L^2}
+||\Pi^\epsilon_0||^2_{L^2}
  +2\sum^{6}_{i=1}R^{\epsilon}_i(t),
\end{align}
where
\begin{align*}
R^{\epsilon}_1(t)=
& {\mu^\epsilon}  \int^t_0\int_{\mathbb{T}^d}\nabla \u^\epsilon\cdot\nabla \u
+       {\nu^\epsilon}  \int^t_0\int_{\mathbb{T}^d}\nabla \H^\epsilon\cdot\nabla \H, \\
R^{\epsilon}_2(t)=&\int_{\mathbb{T}^d}
\big[(\sqrt{\rho^\epsilon}-1)\sqrt{\rho^\epsilon}\u^\epsilon\cdot
 \u \big](t)  -\int_{\mathbb{T}^d}
\big[(\sqrt{\rho^\epsilon}-1)\sqrt{\rho^\epsilon}\u^\epsilon\cdot
 \u \big](0)\nonumber\\
  &     +\int^t_0\int_{\mathbb{T}^d}(\rho^\epsilon-\sqrt{\rho^\epsilon})
  \u^\epsilon\cdot ((\u\cdot \nabla )\u),\\
R^{\epsilon}_3(t)=
& \int^t_0\!\int_{\mathbb{T}^d}(\mathbf{Z}^\epsilon \cdot \nabla)
\H  \cdot [(1-\sqrt{\rho^\epsilon})\u^\epsilon]  - \int^t_0
\!\int_{\mathbb{T}^d}\big\{[(1-\sqrt{\rho^\epsilon})\u^\epsilon
]\cdot \nabla\big\}\H  \cdot   \mathbf{Z}^\epsilon, \nonumber\\
R^{\epsilon}_4(t)=
& \int^t_0\int_{\mathbb{T}^d}(1-\rho^\epsilon)
\u^\epsilon\cdot [(\H\cdot \nabla) \H)]
+\frac12 \int^t_0\int_{\mathbb{T}^d}(\rho^\epsilon-1)\u^\epsilon\cdot \nabla (| \H|^2),\\
R^{\epsilon}_5(t)=&
    \int^t_0 \int_{\mathbb{T}^d}\rho^\epsilon \u^\epsilon \cdot \nabla p,\\
R^{\epsilon}_6(t)=& -\int^t_0\int_{\mathbb{T}^d}(\sqrt{
\rho^\epsilon} \u^\epsilon - \u) \cdot \nabla
\big(\frac{|\u|^2}{2}\big).
\end{align*}

\medskip
\emph{Step 3: Convergence of the  modulated energy functional.}

To show the convergence of the modulated energy functional \eqref{cmi} and to
finish our proof, we have to  estimate the reminders $R^{\epsilon}_i(t)$, $i=1,\dots,6.$

First, in view of \eqref{L2g} and the following two elementary inequalities
\begin{align}
&|\sqrt{x}-1|^2\leq M|x-1|^\gamma,\;\;\;|x-1|\geq\delta, \;\;
\gamma\geq 1,\label{ine10}\\
&|\sqrt{x}-1|^2\leq M|x-1|^2,\;\;\;x\geq 0 \label{ine20}
\end{align}
for some positive constants $M$ and $0<\delta<1$, we obtain
\begin{align}
 \int_{\mathbb{T}^d}|\sqrt{\rho^\epsilon}-1|^2
=&\int_{|\rho^\epsilon-1|\leq
\frac12}|\sqrt{\rho^\epsilon}-1|^2+\int_{|\rho^\epsilon-1|\geq
\frac12}|\sqrt{\rho^\epsilon}-1|^2\nonumber\\
\leq&M\int_{|\rho^\epsilon-1|\leq
\frac12}|\rho^\epsilon-1|^2+M\int_{|\rho^\epsilon-1|\geq
\frac12}|\rho^\epsilon-1|^\gamma\nonumber\\
\leq &M\epsilon^2.\label{ine30}
\end{align}
Now, we begin to estimate the terms $R^\epsilon_i(t),i=1,\dots, 6$.
For the term $R^{\epsilon}_1(t)$, by Young's inequality and the
regularity of $\u$ and $\H$, we have
\begin{align}\label{rr1}
|R^{\epsilon}_1(t)|\leq \frac{\mu^\epsilon}{2} \int^t_0\int_{\mathbb{T}^d}|\nabla \u^\epsilon|^2
+ \frac{\nu^\epsilon}{2}\int^t_0\int_{\mathbb{T}^d}|\nabla \H^\epsilon|^2
+C_T\mu^\epsilon+C_T\nu^\epsilon.
\end{align}

For the  term  $R^{\epsilon}_2(t)$, by H\"older's inequality, the
estimate \eqref{ine30},
  the assumption on the initial data, the estimate on
$\sqrt{\rho^\epsilon}\u^\epsilon$, and the regularity of $\u$, we infer that
\begin{align}\label{rr2}
|R^{\epsilon}_2(t)|\leq  & C\epsilon +||\u(t)||_{L^\infty}
\Big(\int_{\mathbb{T}^d} |\sqrt{\rho^\epsilon}-1|^2 \Big)^{\frac12}
\Big(\int_{\mathbb{T}^d} \rho^\epsilon |\u^\epsilon|^2\Big)^{\frac12}\nonumber\\
&+ ||[(\u\cdot\nabla)\u](t)||_{L^\infty}
\Big(\int^t_0\int_{\mathbb{T}^d} |\sqrt{\rho^\epsilon}-1|^2\Big)^{\frac12}
\Big(\int^t_0\int_{\mathbb{T}^d} \rho^\epsilon |\u^\epsilon|^2 \Big)^{\frac12}\nonumber\\
\leq & C_T\epsilon.
\end{align}

For the term $R^{\epsilon}_3(t)$, making use of the inequality \eqref{ine30},
the basic inequality \eqref{be}, the estimates on $\u^\epsilon$ and
$\H^\epsilon$, the regularity of  $\H$, the assumption \eqref{abc},
and Sobolev's imbedding theorem, we get
\begin{align}\label{rr3}
|R^\epsilon_3(t)|\leq &(||[(\H\cdot \nabla)\H](t)||_{L^\infty}
+ ||\nabla \H(t)||_{L^\infty}\cdot||\H(t)||_{L^\infty})\nonumber\\
& \times\Big(\int^t_0\int_{\mathbb{T}^d}|\sqrt{\rho^\epsilon}-1|^2\Big)^{\frac12}
\Big(\int^t_0\int_{\mathbb{T}^d}|\u^\epsilon|^2\Big)^{\frac12}\nonumber\\
& + ||\nabla \H(t)||_{L^\infty}\int^t_0\Big[\Big(\int_{\mathbb{T}^d}|\sqrt{\rho^\epsilon}-1|^2\Big)^{\frac12}
  ||\u^\epsilon(\tau)||_{L^6}  ||\H^\epsilon(\tau)||_{L^3}\Big]d\tau\nonumber\\
\leq &
C_T\epsilon(1+(\mu^\epsilon)^{-\frac12})+||\sqrt{\rho^\epsilon}-1||_{L^\infty(0,T;L^2)}
\Big(\int^t_0||\u^\epsilon(\tau)||^2_{H^1}d\tau\Big)^{\frac12}
\Big(\int^t_0||\H^\epsilon(\tau)||^2_{H^1}d\tau\Big)^{\frac12}\nonumber\\
\leq & C_T\epsilon(1+(\mu^\epsilon)^{-\frac12}) +C_T\epsilon
(||\u^\epsilon||_{L^2(0,T;L^2)}
+||\nabla\u^\epsilon||_{L^2(0,T;L^2)})\nonumber\\
& \times (||\H^\epsilon||_{L^2(0,T;L^2)}+||\nabla\H^\epsilon||_{L^2(0,T;L^2)})\nonumber\\
\leq & C_T\epsilon(1+(\mu^\epsilon)^{-\frac12})
+C_T\epsilon\big[1+(\mu^\epsilon)^{-\frac12}\big]
\cdot \big[1+(\nu^\epsilon)^{-\frac12}\big]\nonumber\\
\leq & C_T \epsilon^{1-\alpha/2} +C_T \epsilon^{\sigma}\leq C_T
\epsilon^{\sigma},
\end{align}
where $\sigma =1-(\alpha+\beta)/2.$

For the term $R^{\epsilon}_4(t)$, one can utilize the inequality \eqref{ine30},
the estimates on $\u^\epsilon$ and $\sqrt{\rho^\epsilon}\u^\epsilon$, the
regularity of  $\H$, and $\rho^\epsilon-1=\rho^\epsilon-\sqrt{\rho^\epsilon}
+\sqrt{\rho^\epsilon} -1$ to deduce
\begin{align}\label{rr4}
|R^\epsilon_4(t)|\leq &(||[(\H\cdot \nabla)\H](t)||_{L^\infty}+ ||\nabla (|\H|^2)||_{L^\infty})
\Big(\int^t_0\int_{\mathbb{T}^d}|\sqrt{\rho^\epsilon}-1|^2\Big)^{\frac12}\nonumber\\
& \times\bigg[
\Big(\int^t_0\int_{\mathbb{T}^d}|\u^\epsilon|^2\Big)^{\frac12}+
\Big(\int^t_0\int_{\mathbb{T}^d}\rho^\epsilon|\u^\epsilon|^2\Big)^{\frac12} \bigg]\nonumber\\
\leq & C_T \epsilon(1+(\mu^\epsilon)^{-\frac12})\leq  C_T
\epsilon^{1-\alpha/2} .
\end{align}

Using \eqref{ba}, \eqref{L2g} and \eqref{re},
the term $R^{\epsilon}_5(t)$ can be bounded as follows.
\begin{align}\label{rr5}
 | R^{\epsilon}_5(t)|=&
  \Big|\int^t_0 \int_{\mathbb{T}^d}\rho^\epsilon \u^\epsilon \cdot  \nabla p\Big|\nonumber\\
 =&  \Big| \int_{\mathbb{T}^d}\big\{((\rho^\epsilon-1) p)(t)-((\rho^\epsilon-1)
 p)(0)\big\}- \int^t_0\int_{\mathbb{T}^d}(\rho^\epsilon-1)  \partial_t p \Big|\nonumber\\
  \leq &   \Big(\int_{|\rho^\epsilon-1|\leq \frac12}|\rho^\epsilon-1|^2 \Big)^{\frac12}
  \Big[\Big(\int_{\mathbb{T}^3}|p(t)|^2\Big)^{\frac12}+\Big(\int_{\mathbb{T}^3}|p(0)|^2\Big)^{\frac12}\Big]\nonumber\\
  & +  \Big(\int_{|\rho^\epsilon-1|\geq \frac12}|\rho^\epsilon-1|^\gamma \Big)^{\frac{1}{\gamma}}
  \Big[\Big(\int_{\mathbb{T}^3}|p(t)|^{\frac{\gamma}{\gamma-1}}\Big)^{\frac{\gamma-1}{\gamma}}
  +\Big(\int_{\mathbb{T}^3}|p(0)|^{\frac{\gamma}{\gamma-1}}\Big)^{\frac{\gamma-1}{\gamma}}\Big]\nonumber\\
&+\int^t_0 \Big(\int_{|\rho^\epsilon-1|\leq \frac12}|\rho^\epsilon-1|^2\Big)^{\frac12}\Big(\int_{\mathbb{T}^d}|\partial_tp(t)|^2 \Big)^{\frac12}
\nonumber\\
& + \int^t_0 \Big(\int_{|\rho^\epsilon-1|\geq \frac12} |\rho^\epsilon-1|^{\gamma}\Big)^{{\frac{1}{\gamma}}}
\Big(\int_{\mathbb{T}^d}| \partial_t p(t)|^{\frac{\gamma}{\gamma-1}} \Big)^{{\frac{\gamma-1}{\gamma}}}\nonumber\\
\leq & C_T(\epsilon+\epsilon^{2/\kappa})\leq C_T \epsilon,
\end{align}
where $\kappa=\min\{2,\gamma\}$  and we have used the conditions
$s>2+d/2$ and $\gamma>1$.

Finally, to estimate the term  $ R^\epsilon_6(t)$, we rewrite it as
\begin{align}\label{rr6}
  R^{\epsilon}_6(t)=& -\int^t_0\int_{\mathbb{T}^d}(\sqrt{ \rho^\epsilon}\u^\epsilon-
 \u) \cdot \nabla \big(\frac{|\u|^2}{2}\big)\nonumber\\
  =&\int^t_0\int_{\mathbb{T}^d}\sqrt{\rho^\epsilon}(\sqrt{\rho^\epsilon}-1)
 \u^\epsilon \cdot \nabla \big(\frac{|\u|^2}{2}\big)
  -\int^t_0\int_{\mathbb{T}^d}\rho^\epsilon \u^\epsilon \cdot \nabla
 \big(\frac{|\u|^2}{2}\big)\nonumber\\
  =& R^{\epsilon}_{61}(t)+R^{\epsilon}_{62}(t),
\end{align}
where
\begin{align*}
  R^{\epsilon}_{61}(t)& = \int^t_0\int_{\mathbb{T}^d}\sqrt{\rho^\epsilon}(\sqrt{\rho^\epsilon}-1) \u^\epsilon
 \cdot
\nabla \big(\frac{|\u|^2}{2}\big),\\
R^{\epsilon}_{62}(t) &=
   \int^t_0\int_{\mathbb{T}^d}(\rho^\epsilon-1)
 \partial_t\big(\frac{|\u|^2}{2}\big)
   - \int_{\mathbb{T}^d}\Big[\Big((\rho^\epsilon-1)
 \big(\frac{|\u|^2}{2}\big)\Big)(t) -\Big((\rho^\epsilon-1) \big(\frac{|\u|^2}{2}\big)\Big)(0)\Big].
\end{align*}
Applying arguments similar to those used for $R^{\epsilon}_{61}(t)$
and $R^{\epsilon}_{5}(t)$, we arrive at the following boundedness
\begin{align}\label{rr66}
 |R^{\epsilon}_{6}(t)|\leq  |R^{\epsilon}_{61}(t)|+|R^{\epsilon}_{62}(t)|\leq C_T \epsilon.
\end{align}

Inserting the estimates \eqref{rr1}-\eqref{rr66} into \eqref{cmi} and applying Gronwall's inequality,
we conclude
\begin{align}\label{cq}
 &  ||\mathbf{w}^{\epsilon}(t)||^2_{L^2}+
 ||\mathbf{Z}^{\epsilon}(t)||^2_{L^2}+||\Pi^\epsilon(t)||^2_{L^2}\nonumber\\
& \qquad \leq \bar C
\big[||\mathbf{w}^{\epsilon}(0)||^2_{L^2}+||\mathbf{Z}^{\epsilon}(0)||^2_{L^2}
+ ||\Pi^\epsilon_0||^2_{L^2} +C_T\epsilon^\sigma\big],
\quad\mbox{for a.e.}\;\; t\in [0,T],  \end{align} where
\begin{align}\label{ccc} \bar C=\exp {\Big\{C\int^T_0\big[||\nabla
\u(\tau)||_{L^\infty} +
||\nabla\H(\tau)||_{L^\infty}\big]d\tau\Big\}}<+\infty.
\end{align}

Now, letting $\epsilon $ go to $0$, we obtain $\mathbf{K}=\H$  in
$L^\infty(0,T; L^2)$ and $\mathbf{J}=\u$ in $L^\infty(0,T; L^2)$.
The inequality \eqref{icd20} follows from \eqref{icd2} and \eqref{cq} directly.
Thus, we complete the proof.
\end{proof}

\begin{proof}[Proof of Theorem \ref{MRb}]
 For simplicity we assume here that $\mu^\epsilon\equiv\mu$,
$\lambda^\epsilon\equiv\lambda$, and $\nu^\epsilon\equiv\nu$ are
constants, independent of $\epsilon$, satisfying  $\mu>0$,
$\mu+\lambda>0$, and $\nu>0$. The case \eqref{pa} can be treated
similarly. The proof of Theorem \ref{MRb} is similar to that of
Theorem \ref{MRa}. Since the viscosity is involved here, we have to
modulate the part of dissipation energy in the energy inequality
\eqref{be}. We state the main different points in the proof here.

From the basic energy inequality \eqref{be}, we obtain that, for a.e.
$t\in [0,T]$, $\rho^\epsilon |\u^\epsilon|^2$ and
$\big((\rho^\epsilon)^\gamma-1-\gamma(\rho^\epsilon-1)\big)/{\epsilon^2}$
are bounded in $L^\infty(0,T;L^1)$, $ \H^\epsilon $ is
bounded in $L^\infty(0,T;L^2)$,    $\nabla \u^\epsilon $ is
bounded in $L^2(0,T;L^2)$,  and $\nabla \H^\epsilon $ is
bounded in $L^2(0,T;L^2)$. Therefore, we have
\begin{equation*}
    \rho^\epsilon \rightarrow    1 \ \
    \text{strongly in}\ \
C([0,T],L^\gamma_2(\mathbb{T}^d)),
\end{equation*}
and $\u^\epsilon$ is bounded in $L^2(0,T;L^2)$.
The boundedness of $\rho^\epsilon |\u^\epsilon|^2$ and $ |\H^\epsilon|^2$
in $L^\infty(0,T;L^1)$ implies the following convergence (up to
the extraction of a subsequence $\epsilon_n$):
\begin{align*}
& \sqrt{\rho^\epsilon}  \u^\epsilon   \  \text{converges weakly-$\ast$}
\,\, \text{to some}\, \bar{\mathbf{J}} \, \,\text{in}\,
L^\infty(0,T;L^2(\mathbb{T}^d)),\\
&  \H^\epsilon \    \text{converges weakly-$\ast$}
\,\, \text{to some}\,  \bar{\mathbf{K}} \, \,\text{in}\,
L^\infty(0,T;L^2(\mathbb{T}^d)).
\end{align*}
Our main task in this section is to show that
$ \bar{\mathbf{J}}=\u$ and  $ \bar{\mathbf{K}}=\H$ in some
sense, where $(\u,\H)$  is the strong solution to the
viscous incompressible  MHD equations \eqref{a2ll}-\eqref{a2nn}.

  Next, we shall also  modulate the
energy inequality \eqref{be}. The conservation of energy for the
viscous incompressible  MHD equations \eqref{a2ll}-\eqref{a2nn} reads
\begin{equation}\label{cff1}
  \frac12 \int_{\mathbb{T}^d}\big[|\u|^2 + |\H|^2\big](t)
  +\int^t_0 \int_{\mathbb{T}^d} \big[\mu|\nabla\u|^2+ \nu|\nabla\H|^2\big]
   =\frac12  \int_{\mathbb{T}^d}  \big[|\u_0|^2+|\H_0|^2\big].
\end{equation}
Similarly to Step 2, we use $\u$ to test the momentum equation \eqref{a2j} to deduce
\begin{align}\label{civ1}
&  \int_{\mathbb{T}^d}(\rho^\epsilon \u^\epsilon \cdot \u)(t)
+\int^t_0\int_{\mathbb{T}^d}\rho^\epsilon \u^\epsilon \cdot \big[(\u\cdot \nabla)
\u-(\H\cdot \nabla)\H-\mu \Delta \u+\nabla p+\frac12\nabla(|\H|^2)\big] \nonumber\\
&  -\int^t_0\int_{\mathbb{T}^d}\big[(\rho^\epsilon \u^\epsilon \otimes
\u^\epsilon)\cdot \nabla \u+   (\H^\epsilon\cdot \nabla)\H^\epsilon \cdot \u
 - \mu   \nabla  \u^\epsilon  \cdot \nabla \u\big] = \int_{\mathbb{T}^d}
  \rho^\epsilon_0 \u^\epsilon_0 \cdot \u_0.
\end{align}
Then, we test \eqref{a2k} by $\H$ to infer that
\begin{align}\label{cjj1}
&  \int_{\mathbb{T}^d}(\H^\epsilon \cdot \H)(t)
+\int^t_0\int_{\mathbb{T}^d}\H^\epsilon \cdot\big[(\u\cdot \nabla)
\H-(\H\cdot \nabla)\u-\nu \Delta \H\big] + \nu
\int^t_0\int_{\mathbb{T}^d}
\nabla \H^\epsilon  \cdot \nabla \H\nonumber\\
& \quad +\int^t_0\int_{\mathbb{T}^d}\big[(\dv  {\u}^\epsilon) { \H}^\epsilon +
( {\u}^\epsilon \cdot \nabla) {\H}^\epsilon  - ( {\H}^\epsilon\cdot \nabla)
 {\u}^\epsilon\big]\cdot \H  = \int_{\mathbb{T}^d}\H^\epsilon_0 \cdot \H_0.
\end{align}
Summing up \eqref{cei} and   \eqref{cff1}, and inserting
\eqref{cii}, \eqref{cj1i} with $\mu^\epsilon\equiv\mu$ and
$\lambda^\epsilon\equiv\lambda$, \eqref{civ1}, and \eqref{cjj1} into
the resulting inequality, we deduce the following inequality by a
straightforward calculation
\begin{align}\label{ckk}
& \ \  \frac{1}{2}\int_{\mathbb{T}^d}\Big\{  |\sqrt{\rho^\epsilon} \u^\epsilon
-\u|^2(t)+ |  \H^\epsilon
-\H|^2(t)
  +(\Pi^\epsilon)^2(t)\Big\}\nonumber\\
&\quad
+ \mu\int^t_0\int_{\mathbb{T}^d}|\nabla \u^\epsilon -\nabla \u |^2
+ \frac{\nu}{2}  \int^t_0\int_{\mathbb{T}^d}|\nabla \H^\epsilon-\nabla \H |^2
+ (\mu+\lambda)\int^t_0\! \int_{\mathbb{T}^d}|{\rm div}\u^\epsilon |^2\nonumber\\
&\quad
+ \frac{\nu}{2}  \int^t_0\int_{\mathbb{T}^d}|\nabla \H^\epsilon|^2
+ \frac{\nu}{2}  \int^t_0\int_{\mathbb{T}^d}|\nabla \H|^2
\nonumber\\
&\leq  -\int^t_0\int_{\mathbb{T}^d}\rho^\epsilon
\u^\epsilon\cdot [(\H\cdot \nabla) \H)] -\int^t_0\int_{\mathbb{T}^d}
(\H^\epsilon\cdot \nabla) \H^\epsilon \cdot\u           \nonumber\\
&\quad       +\int^t_0\int_{\mathbb{T}^d} \H^\epsilon \cdot\big[(\u\cdot \nabla)
\H-(\H\cdot \nabla)\H\big]  +\frac12\int^t_0\int_{\mathbb{T}^d}
\rho^\epsilon \u^\epsilon \cdot \nabla(|\H|^2)\nonumber\\
&\quad  +\int^t_0\int_{\mathbb{T}^d}\big[(\dv  {\u}^\epsilon) { \H}^\epsilon
 + ( {\u}^\epsilon \cdot \nabla) {\H}^\epsilon
 - ( {\H}^\epsilon\cdot \nabla) {\u}^\epsilon\big]\cdot \H
+\mu\int^t_0\int_{\mathbb{T}^d}(1-\rho^\epsilon) \u^\epsilon \Delta u
  \nonumber\\
&\quad +\int^t_0\int_{\mathbb{T}^d}\rho^\epsilon \u^\epsilon \cdot\big[ (\u\cdot
\nabla) \u +\nabla p \big]-\int^t_0\int_{\mathbb{T}^d}(\rho^\epsilon
\u^\epsilon\otimes \u^\epsilon)\cdot \nabla \u\nonumber\\
&\quad
+\int_{\mathbb{T}^d}
\big[(\sqrt{\rho^\epsilon}-1)\sqrt{\rho^\epsilon}\u^\epsilon\cdot
\u\big](t)
    -\int_{\mathbb{T}^d}
\big[(\sqrt{\rho^\epsilon}-1)\sqrt{\rho^\epsilon}\u^\epsilon\cdot
\u\big](0)\nonumber\\
&\quad   +\frac{1}{2}\int_{\mathbb{T}^d}\Big\{  |\sqrt{\rho^\epsilon} \u^\epsilon
-\u|^2(0)+|\H^\epsilon
-\H |^2(0)
  +(\Pi^\epsilon_0)^2\Big\}.
\end{align}

By virtue of \eqref{ck1i} and \eqref{cli}, we can rewrite the inequality \eqref{ckk}
as follows
\begin{align}\label{ckk1}
& \ \  \frac{1}{2}\int_{\mathbb{T}^d}\Big\{  |\sqrt{\rho^\epsilon} \u^\epsilon
-\u|^2(t)+ |  \H^\epsilon
-\H|^2(t)
  +(\Pi^\epsilon)^2(t)\Big\}\nonumber\\
&\quad
+  {\mu}   \int^t_0\int_{\mathbb{T}^d}|\nabla \u^\epsilon -\nabla \u |^2
+ \frac{\nu}{2}  \int^t_0\int_{\mathbb{T}^d}|\nabla \H^\epsilon-\nabla \H |^2\nonumber\\
&\quad + (\mu+\lambda)\int^t_0\! \int_{\mathbb{T}^d}|{\rm div}\u^\epsilon|^2
+ \frac{\nu}{2}  \int^t_0\int_{\mathbb{T}^d}|\nabla \H^\epsilon|^2
+ \frac{\nu}{2}  \int^t_0\int_{\mathbb{T}^d}|\nabla \H|^2
\nonumber\\
&\leq
     \frac{1}{2}\int_{\mathbb{T}^d}\Big\{  |\sqrt{\rho^\epsilon} \u^\epsilon
-\u|^2(0)+|\H^\epsilon
-\H |^2(0)
  +(\Pi^\epsilon_0)^2\Big\}\nonumber\\
  & \quad + R^\epsilon_2(t)+ R^\epsilon_4(t)+ R^\epsilon_5(t)+ R^\epsilon_6(t)+ R^\epsilon_7(t)+ R^\epsilon_8(t),
\end{align}
where $R^\epsilon_2(t)$ and $R^\epsilon_i(t)$, $i=4,5,6$, are the same as before, and
\begin{align*}
     R^\epsilon_7(t)&=\mu\int^t_0\int_{\mathbb{T}^d}(1-\rho^\epsilon) \u^\epsilon \Delta u,\\
    R^\epsilon_8(t) &=
  \int^t_0\!\int_{\mathbb{T}^d}(\mathbf{Z}^\epsilon \cdot \nabla)
\H  \cdot [(1-\sqrt{\rho^\epsilon})\u^\epsilon]  - \int^t_0\!\int_{\mathbb{T}^d}\big\{[(1-\sqrt{\rho^\epsilon})\u^\epsilon
]\cdot \nabla\big\}\H  \cdot   \mathbf{Z}^\epsilon.
\end{align*}

Form the previous arguments on $R^\epsilon_2(t)$ and $R^\epsilon_i(t)$, $i=4,5,6$, we get
\begin{align}\label{ckk2}
 |R^\epsilon_2(t)|+\sum^6_{i=4}| R^\epsilon_i(t)|\leq C_T \epsilon.
\end{align}

Now, we estimate $R^\epsilon_7(t)$ and $R^\epsilon_8(t)$. Using  the inequality \eqref{ine30},
H\"older's inequality, the estimates on $\u^\epsilon$ and
$\sqrt{\rho^\epsilon}\u^\epsilon$,
  the
regularity of  $\u$, and $\rho^\epsilon-1=\rho^\epsilon-\sqrt{\rho^\epsilon}+\sqrt{\rho^\epsilon} -1$, we obtain
\begin{align}\label{ckk3}
|R^\epsilon_7(t)|\leq & \mu ||\Delta\u(t)||_{L^\infty}
\Big(\int^t_0\int_{\mathbb{T}^d}|\sqrt{\rho^\epsilon}-1|^2\Big)^{\frac12}
 \bigg[
\Big(\int^t_0\int_{\mathbb{T}^d}|\u^\epsilon|^2\Big)^{\frac12}+
\Big(\int^t_0\int_{\mathbb{T}^d}\rho^\epsilon|\u^\epsilon|^2\Big)^{\frac12} \bigg]\nonumber\\
\leq & C_T \epsilon(1+\frac{1}{\sqrt{\mu}})\leq C_T
\frac{\epsilon}{\sqrt{\mu}}.
\end{align}

For the term $R^{\epsilon}_8(t)$, we can make use of \eqref{ine30}, \eqref{be}, and
the estimates on $\u^\epsilon$ and $\H^\epsilon$, the regularity of $\H$,
the assumption \eqref{abc}, and Sobolev's imbedding theorem to deduce
\begin{align}\label{ckk4}
|R^\epsilon_8(t)|\leq &(||[(\H\cdot \nabla)\H](t)||_{L^\infty}
+ ||\nabla \H(t)||_{L^\infty}\cdot||\H(t)||_{L^\infty})\nonumber\\
& \times\Big(\int^t_0\int_{\mathbb{T}^d}|\sqrt{\rho^\epsilon}-1|^2\Big)^{\frac12}
\Big(\int^t_0\int_{\mathbb{T}^d}|\u^\epsilon|^2\Big)^{\frac12}\nonumber\\
& + ||\nabla \H(t)||_{L^\infty}\int^t_0\Big[\Big(\int_{\mathbb{T}^d}|
\sqrt{\rho^\epsilon}-1|^2\Big)^{\frac12}
  ||\u^\epsilon(\tau)||_{L^6}  ||\H^\epsilon(\tau)||_{L^3}\Big]d\tau\nonumber\\
\leq &
C_T\epsilon(1+\frac{1}{\sqrt{\mu}})+||\sqrt{\rho^\epsilon-1}||_{L^\infty(0,T;L^2)}
\Big(\int^t_0||\u^\epsilon(\tau)||^2_{H^1}d\tau\Big)^{\frac12}
\Big(\int^t_0||\H^\epsilon||^2_{H^1}(\tau)d\tau\Big)^{\frac12}\nonumber\\
\leq & C_T\epsilon(1+\frac{1}{\sqrt{\mu}}) +C_T\epsilon
(||\u^\epsilon||_{L^2(0,T;L^2)}
+||\nabla\u^\epsilon||_{L^2(0,T;L^2)})\nonumber\\
& \times
(||\H^\epsilon||_{L^2(0,T;L^2)}+||\nabla\H^\epsilon||_{L^2(0,T;L^2)})\nonumber\\
\leq & C_T\epsilon(1+\frac{1}{\sqrt{\mu}}) +C_T\epsilon\big[1+\mu^{-\frac12}\big] \cdot \big[1+\nu^{-\frac12}\big]\nonumber\\
\leq & C_T \epsilon +C_T \epsilon/\sqrt{\mu\nu}\leq C_T
\epsilon/\sqrt{\mu\nu}.
\end{align}

Now, substituting \eqref{ckk2}-\eqref{ckk4} into \eqref{ckk1}
and applying Gronwall's inequality, we conclude
\begin{align}\label{ckk5}
 &  ||\mathbf{w}^{\epsilon}(t)||^2_{L^2}+
 ||\mathbf{Z}^{\epsilon}(t)||^2_{L^2}+||\Pi^\epsilon(t)||^2_{L^2}\nonumber\\
& \quad \leq \bar C
\big[||\mathbf{w}^{\epsilon}(0)||^2_{L^2}+||\mathbf{Z}^{\epsilon}(0)||^2_{L^2}
+ ||\Pi^\epsilon_0||^2_{L^2} +C_T\epsilon/\sqrt{\mu\nu}\big],
\quad\mbox{for a.e.}\;\; t\in [0,T],
\end{align}
where $\bar C$ is defined by \eqref{ccc}.  Combining \eqref{icd22}
with \eqref{ckk5} we obtain \eqref {icd220}. Substituting \eqref{ckk5} into \eqref{ckk1},
we conclude that $\nabla \u^\epsilon$ converges to $\nabla \u$ strongly
 in $L^2 (0,T;L^2(\mathbb{T}^d))$
 and  $\nabla \H^\epsilon$ to $\nabla \H$ strongly
 in $L^2 (0,T;  L^2(\mathbb{T}^d))$. This completes the proof of Theorem \ref{MRb}.
\end{proof}

\section{Proof of  Theorem  \ref{MRc}}

In this section we shall study the incompressible limit of the
compressible MHD equations \eqref{a2i}-\eqref{a2k} with general
initial data. Compared with the case of the well-prepared initial
data, the main difficulty here is to control the oscillations caused
by the initial data. For simplicity, we assume here that
$\mu^\epsilon\equiv\mu$, $\lambda^\epsilon\equiv\lambda$, and
$\nu^\epsilon\equiv\nu$ are constants, independent of $\epsilon$,
satisfying  $\mu>0$,  $2\mu+d\lambda>0$, and $\nu>0$.

\begin{proof}[Proof of Theorem  \ref{MRc}]

As stated in the proof of Theorem \ref{MRb},
  we obtain from the basic energy inequality \eqref{be} that, for a.e.
  $t\in [0,T]$, $\rho^\epsilon |\u^\epsilon|^2$ and
$ \big((\rho^\epsilon)^\gamma-1-\gamma(\rho^\epsilon-1)\big)/{\epsilon^2}$
are bounded in $L^\infty(0,T;L^1)$, $ \H^\epsilon $ is
bounded in $L^\infty(0,T;L^2)$,    $\nabla \u^\epsilon $ is
bounded in $L^2(0,T;L^2)$,  and $\nabla \H^\epsilon $ is
bounded in $L^2(0,T;L^2)$. Therefore, we have
\begin{equation}\label{rhoc}
    \rho^\epsilon \rightarrow    1 \ \
    \text{strongly in}\ \
C([0,T],L^\gamma_2(\mathbb{T}^d)),
\end{equation}
and $\u^\epsilon$ is bounded in $L^2(0,T;L^2)$.
The fact that $\rho^\epsilon |\u^\epsilon|^2$ and $ |\H^\epsilon|^2$ are bounded
in $L^\infty(0,T;L^1)$ gives the following convergence (up to
the extraction of a subsequence $\epsilon_n$):
\begin{align*}
& \sqrt{\rho^\epsilon}  \u^\epsilon   \  \text{converges weakly-$\ast$}
\,\, \text{to some}\,  { \bar{\mathbf{J}}} \, \,\text{in}\,
L^\infty(0,T;L^2(\mathbb{T}^d)),\\
&  \H^\epsilon \    \text{converges weakly-$\ast$}
\,\, \text{to some}\,  { \bar{\mathbf{K}}} \, \,\text{in}\,
L^\infty(0,T;L^2(\mathbb{T}^d)).
\end{align*}
Our main task in this section is to show that
$    \bar{\mathbf{J}} =\u$ and  $    \bar{\mathbf{K}}  =\H$ in some
sense, where $(\u,\H)$  is the strong solution to the
incompressible viscous MHD equations \eqref{a2ll}-\eqref{a2nn}.
The key point is to control the oscillations caused by the
  initial data. This can be done as follows.

\medskip
\emph{Step 1: Description and cancelation of the oscillations.}

In order to describe the oscillations caused by the initial data,
we employ the ``filtering" method which has been used previously by several authors,
see~\cite{LM98,G,M01b,D02}.

We project the momentum equation  \eqref{a2j} on the
``gradient vector-fields'' to  find
\begin{align}\label{Pr}
&\partial_t Q(\rho^\epsilon\u^\epsilon)+Q[\dv(\rho^\epsilon\u^\epsilon\otimes\u^\epsilon)]
  - (2\mu+\lambda)\nabla\dv\u^\epsilon+\frac12\nabla(|{\H}^\epsilon|^2)\nonumber\\
 &\quad -Q[({\H}^\epsilon\cdot \nabla
  ){\H}^\epsilon] +\frac{a}{\epsilon^2}\nabla \big((\rho^\epsilon)^\gamma
 -1-\gamma (\rho^\epsilon-1)\big)+\frac{1}{\epsilon^2}\nabla (\rho^\epsilon-1)=0.
\end{align}
Noticing $\rho^\epsilon=1+\epsilon\varphi^\epsilon $,  we can write   \eqref{a2i} and \eqref{Pr} as
\begin{align}
&\epsilon\partial_t \varphi^\epsilon+\dv Q(\rho^\epsilon \u^\epsilon)=0,\label{ep1}\\
 &\epsilon\partial_t  Q(\rho^\epsilon \u^\epsilon)+ \nabla \varphi^\epsilon=\epsilon\mathbf{F}^\epsilon,\label{ep2}
\end{align}
where $ \mathbf{F}^\epsilon$ is given by
\begin{align}
\mathbf{F}^\epsilon=&-Q[\dv(\rho^\epsilon\u^\epsilon\otimes\u^\epsilon)]
  +(2\mu+\lambda)\nabla\dv\u^\epsilon-\frac12\nabla(|{\H}^\epsilon|^2)\nonumber\\
  & +Q[({\H}^\epsilon\cdot \nabla
  ){\H}^\epsilon] -\frac{a}{\epsilon^2}\nabla \big((\rho^\epsilon)^\gamma
 -1-\gamma (\rho^\epsilon-1)\big). \label{ff}
\end{align}
Therefore, we introduce the following group defined by
$\mathcal{L}(\tau)=e^{\tau L}$, $\tau \in \mathbb{R}$, where $L$ is
the operator defined on $\mathcal{D}_0'\times (\mathcal{D}')^{d}$
with $\mathcal{D}_0'=\{\phi \in \mathcal{D}',
\int_{\mathbb{T}^d}\phi(x)dx=0\}$, by
$$
L\Big(\begin{array}{c}
 \phi   \\
\mathbf{v}
\end{array}
\Big) =\Big(\begin{array}{c}
 -\text{div}\,\mathbf{v}  \\
-\nabla \phi
\end{array}
\Big).
$$
Then, it is easy to check that $e^{\tau L}$ is an isometry on each
 $H^r\times (H^r)^d$ for all $r\in \mathbb{R}$ and for all $\tau \in
 \mathbb{R}$.
Denoting
$$
\Big(\begin{array}{c}
 \bar \phi (\tau)   \\
\bar{\mathbf{v}}(\tau)
\end{array}
\Big) =e^{\tau L}\Big(\begin{array}{c}
 \phi  \\
\mathbf{v}
\end{array}
\Big),
$$
we have
$$
\frac{\partial \bar \phi}{\partial \tau}=-\text{div} \bar{\mathbf{v}}, \quad
 \frac{\partial \bar{\mathbf{v}}}{\partial \tau}=-\nabla   \bar \phi.
$$
Thus, $\frac{\partial^2\bar \phi}{\partial \tau ^2}-\Delta \bar \phi
 =0$.

 In the sequel, we shall denote
 $$\mathbf{U}^\epsilon =\left(\begin{array}{c}
                           \varphi^\epsilon \\
                          Q(\rho^\epsilon \u^\epsilon)
                        \end{array}\right),\quad
\mathbf{V}^\epsilon= \mathcal{L}\Big(-\frac{t}{\epsilon}\Big)
\left(\begin{array}{c}
 \varphi^\epsilon\\
 Q(\rho^\epsilon \u^\epsilon)
 \end{array}\right)
 $$
and use the following approximations
  $$\bar{\mathbf{U}}^\epsilon =\left(\begin{array}{c}
                           \Phi^\epsilon \\
                          Q(\sqrt{\rho^\epsilon} \u^\epsilon)
                        \end{array}\right),\quad
\bar{\mathbf{V}}^\epsilon= \mathcal{L}\Big(-\frac{t}{\epsilon}\Big)
\left(\begin{array}{c}
 \Phi^\epsilon\\
 Q(\sqrt{\rho^\epsilon} \u^\epsilon)
 \end{array}\right),
 $$
which satisfy
\begin{align}\label{UU}
 ||{\mathbf{U}}^\epsilon-\bar{\mathbf{U}}^\epsilon||_{L^\infty(0,T;
 L^{\frac{2\gamma}{\gamma+1}}(\mathbb{T}^d))}\rightarrow 0  \ \
  \text{as}  \ \ \epsilon \rightarrow 0.
\end{align}

With this notation, we can rewrite the equations \eqref{ep1}-\eqref{ep2} as
$$
\partial_t \mathbf{U}^\epsilon =\frac{1}{\epsilon} L \mathbf{U}^\epsilon +\widehat{\mathbf{F}^\epsilon},
$$
or equivalently
\begin{align}\label{ep3}
\partial_t \mathbf{V}^\epsilon =  \mathcal{L}\Big(-\frac{t}{\epsilon}\Big) \widehat{\mathbf{F}^\epsilon},
\end{align}
where (and in what follows) $\widehat{\mathbf{v}}  $ denotes $(0, \mathbf{v})^\mathrm{T}$.

It is easy to check that $\mathbf{F}^\epsilon$, given by \eqref{ff}, is bounded
in $L^2(0,T;H^{-s_0}(\mathbb{T}^d))$ for some $s_0$ ($ s_0\in \mathbb{R})$. Hence,
$\mathbf{V}^\epsilon$ is compact in time. Moreover, by virtue of the energy inequality \eqref{be}
and the boundedness  of the linear projector $P$,
$\mathbf{V}^\epsilon \in L^\infty(0,T; L^{\frac{2\gamma}{\gamma+1}}(\mathbb{T}^d))$
uniformly in $\epsilon$. Thus,
\begin{align}\label{ep4}
\mathbf{V}^\epsilon \,\,  \text{converges strongly to some}\,
 \,  {\bar{\mathbf{V}}}   \,\, \text{in}\,\,
L^r(0,T;H^{-s'}(\mathbb{T}^d))
\end{align}
 for all $ s'> s_0$ and $1 < r<\infty$.

Denote $\theta\equiv 2\mu+\lambda$, $\mathcal{L}_1(\tau)$ the first
component of $ \mathcal{L}(\tau)$, and $\mathcal{L}_2(\tau)$ the
last $d$ components of $ \mathcal{L}(\tau)$. If we had sufficient
compactness in space, then we could pass the limit in \eqref{ep3}
and obtain the following limit system for the oscillating parts
\begin{align}\label{osci}
  \partial_t \bar{\mathbf{V}} +\mathcal{Q}_1(\u, \bar{\mathbf{V}})
  +\mathcal{Q}_2(\bar{\mathbf{V}}, \bar{\mathbf{V}})-\frac{\theta}{2}\Delta   \bar{\mathbf{V}}=0,
\end{align}
where $\u$ is the strong solution of the viscous incompressible MHD
 equations \eqref{a2ll}-\eqref{a2nn}, $\mathcal{Q}_1$ is a linear
form of ${\mathbf{V}}$ defined by
\begin{align}\label{ep5}
\mathcal{Q}_1(\v, \mathbf{V})=\lim_{\tau \rightarrow \infty}\frac{1}{\tau} \int^\tau_0 \mathcal{L}(-s)
\bigg(
\begin{array}{c}
   0 \\
  \dv(\v\otimes \mathcal{L}_2(s)\mathbf{V}+\mathcal{L}_2(s)\mathbf{V}\otimes\v)
\end{array}
\bigg)ds,
\end{align}
and  $\mathcal{Q}_2$ is  a bilinear form of ${\mathbf{V}}$ defined by
\begin{align}\label{ep6}
\mathcal{Q}_2(\mathbf{V}, \mathbf{V})=\lim_{\tau \rightarrow \infty}\frac{1}{\tau} \int^\tau_0 \mathcal{L}(-s)
\bigg(
\begin{array}{c}
   0 \\
  \dv( \mathcal{L}_2(s)\mathbf{V}\otimes \mathcal{L}_2(s)\mathbf{V})
  +\frac{\gamma-1}{2}\nabla (\mathcal{L}_1(s)\mathbf{V})^2
\end{array}
\bigg)ds
\end{align}
for any divergence-free vector field $\v \in L^2(\mathbb{T}^d)^d$
and any $\mathbf{V} =(\phi, \nabla q)^{\mathrm{T}} \in
L^2(\mathbb{T}^{d})^{d+1}$. Actually, the convergence in \eqref{ep5} and
\eqref{ep6} can be guaranteed by the following Proposition.
\begin{prop}[\cite{M01b}]\label{P2}
For all $\v\in L^{r_1}(0,T;L^2)$ and $\mathbf{V}\in L^{r_2}(0,T;L^2)$, we have the
following weak convergences (${r_1}$ and ${r_2}$ are such that the products are well defined)
\begin{align}
 & w-\lim_{\epsilon\rightarrow 0} \mathcal{L}
 \Big(-\frac{t}{\epsilon}\Big)\bigg(
\begin{array}{c}
   0 \\
  \dv(\v\otimes \mathcal{L}_2(\frac{t}{\epsilon})\mathbf{V}
  +\mathcal{L}_2(\frac{t}{\epsilon})\mathbf{V}\otimes\v)
\end{array}
\bigg)=  \mathcal{Q}_1(\v,\mathbf{V}),  \label{ep13}\\
 & w-\lim_{\epsilon\rightarrow 0} \mathcal{L}
 \Big(-\frac{t}{\epsilon}\Big)\bigg(
 \begin{array}{c}
   0 \\
  \dv( \mathcal{L}_2(\frac{t}{\epsilon})\mathbf{V}\otimes \mathcal{L}_2(\frac{t}{\epsilon})\mathbf{V})
  +\frac{\gamma-1}{2}\nabla (\mathcal{L}_1(\frac{t}{\epsilon})\mathbf{V})^2
\end{array}
\bigg) = \mathcal{Q}_2(\mathbf{V},\mathbf{V}). \label{ep14}
\end{align}
\end{prop}
The viscosity term in the oscillation equations \eqref{osci} is obtained by the following Proposition.

\begin{prop}[\cite{M01b}]\label{Pv}
Suppose that the same hypothesis as in Proposition \ref{P2} on $\mathbf{V}$ holds. Then, we have
\begin{align}\label{ep7}
 \frac{\theta}{2}\Delta \mathbf{V}=\lim_{\tau \rightarrow \infty}\frac{1}{\tau} \int^\tau_0 \mathcal{L}(-s)
\bigg(
\begin{array}{c}
   0 \\
  \theta \Delta  \mathcal{L}_2(s)\mathbf{V}
\end{array}
\bigg)ds.
\end{align}
\end{prop}

The following propositions, the proof of which can be found in
\cite{M01b}, play an important role in our subsequent analysis.

\begin{prop}[\cite{M01b}]\label{P1}
For all $\v,\mathbf{V}, \mathbf{V}_1 $ and $ \mathbf{V}_2 $ (regular enough to define
all the products), we have
\begin{align}
 &   \int_{\mathbb{T}^d} \mathcal{Q}_1(\v,\mathbf{V} )\mathbf{V} =0, \ \
   \int_{\mathbb{T}^d} \mathcal{Q}_2(\mathbf{V},\mathbf{V} )\mathbf{V} =0, \label{ep10}\\
& \int_{\mathbb{T}^d} [\mathcal{Q}_1(\v,\mathbf{V}_1)\mathbf{V}_2+
   \mathcal{Q}_1(\v,\mathbf{V}_2)\mathbf{V}_1]=0, \label{ep11}\\
 & \int_{\mathbb{T}^d} [\mathcal{Q}_2(\mathbf{V}_1,\mathbf{V}_1)\mathbf{V}_2+
   2\mathcal{Q}_2(\mathbf{V}_1,\mathbf{V}_2)\mathbf{V}_1]=0. \label{ep12}
\end{align}
\end{prop}

\begin{prop}[\cite{M01b}]\label{P3}
Using the symmetry of $\mathcal{Q}_2$, we can extend the equality
\eqref{ep14} in Proposition \ref{P2}  to the case:
\begin{align}\label{ep77}
 & w-\lim_{\epsilon\rightarrow 0} \mathcal{L}
 \Big(-\frac{t}{\epsilon}\Big)\bigg\{\bigg(
\begin{array}{c}
   0 \\
  \frac12\dv\big[\mathcal{L}_2(\frac{t}{\epsilon})\mathbf{V}_1\otimes
  \mathcal{L}_2(\frac{t}{\epsilon})\mathbf{V}_2 +\mathcal{L}_2(\frac{t}{\epsilon})
  \mathbf{V}_2\otimes \mathcal{L}_2(\frac{t}{\epsilon})\mathbf{V}_1\big]
\end{array}
\bigg)   \nonumber\\
  &\qquad \qquad \qquad \quad \ \ +  \bigg(
\begin{array}{c}
   0 \\
   \frac{\gamma-1}{2}\nabla (\mathcal{L}_1(\frac{t}{\epsilon})\mathbf{V}_1 \otimes
   \mathcal{L}_1(\frac{t}{\epsilon})\mathbf{V}_2)
  \end{array}\bigg)\bigg\}
   =  \mathcal{Q}_2(\mathbf{V}_1,\mathbf{V}_2).
\end{align}
Moreover, the above identity holds for $\mathbf{V}_1\in L^q(0,T;H^r)$ and  $\mathbf{V}_2\in L^p(0,T;H^{-r})$
with $r\in \mathbb{R}$ and $1/p +1/q=1$. Also, (\ref{ep77}) can be extended to the case where we replace
$\mathbf{V}_2$ in the left-hand side  by a sequence $\mathbf{V}^\epsilon_2$ such that $\mathbf{V}^\epsilon_2$
converges strongly to $\mathbf{V}_2$ in  $ L^p(0,T;H^{-r})$.
\end{prop}

\medskip
\emph{Step 2: The modulated energy functional and uniform estimates.}

 Let $\mathbf{V}^0$ be the solution of the following system
 \begin{align}\label{ep8}
 \partial_t  \mathbf{V}^0 +\mathcal{Q}_1(\u, \mathbf{V}^0)
  +\mathcal{Q}_2(\mathbf{V}^0, \mathbf{V}^0)-\frac{\theta}{2}\Delta  \mathbf{V}^0=0
 \end{align}
 with initial data
 \begin{align}\label{ep9}
{\mathbf{V}^0}|_{t=0}=(\varphi_0,Q\tilde{\u}_0)^{\mathrm{T}},
 \end{align}
where $\u$ is the strong solution of  the viscous incompressible MHD equations \eqref{a2ll}-\eqref{a2nn} with initial velocity $\u_0$.
From \cite{M01b}, we know that the Cauchy problem \eqref{ep8}-\eqref{ep9} has a unique global
strong solution.

In order to prove the convergence results in Theorem \ref{MRc}, we have to bound the term
$$
\Big\|\sqrt{\rho^\epsilon} \u^\epsilon -\u-\mathcal{L}_2\Big(\frac{t}{\epsilon}\Big)\mathbf{V}^0\Big\|^2_{L^2{(\mathbb{T}^d})}
+\|\H^\epsilon-\H\|^2_{L^2{(\mathbb{T}^d})}
+\Big\|\Pi^\epsilon-\mathcal{L}_1\Big(\frac{t}{\epsilon}\Big)\mathbf{V}^0\Big\|^2_{L^2({\mathbb{T}^d})}.
$$

To this end, we first recall the following energy inequality of the compressible
MHD equations \eqref{a2i}-\eqref{a2k}:
 \begin{align}\label{ce}
 & \frac12\int_{\mathbb{T}^d}\Big[\rho^\epsilon(t)|\u^\epsilon|^2(t)+|\H^\epsilon|^2(t)
    +(\Pi^\epsilon(t))^2
      \Big]
+\mu \int^t_0\! \int_{\mathbb{T}^d}|\nabla \u^\epsilon|^2 \nonumber\\
&\ \  +(\mu+\lambda)\int^t_0\! \int_{\mathbb{T}^d}|{\rm div}\u^\epsilon |^2
+\nu\int^t_0\! \int_{\mathbb{T}^d}|\nabla \H^\epsilon |^2\nonumber\\
& \ \ \   \leq
\frac12\int_{\mathbb{T}^d}\Big[\rho_0^\epsilon|\u_0^\epsilon|^2
+|\H_0^\epsilon|^2+(\Pi^\epsilon_0)^2     \Big], \quad\mbox{for
a.e.}\;\; t\in [0,T].
 \end{align}
On the other hand, the conservation of energy for the
incompressible viscous MHD equations \eqref{a2ll}-\eqref{a2nn}
reads
\begin{equation}\label{cff}
  \frac12 \int_{\mathbb{T}^d}\big[|\u|^2(t)+ |\H|^2(t)\big]
  +\int^t_0 \int_{\mathbb{T}^d} \big[\mu|\nabla\u |^2+ \nu|\nabla\H |^2\big]
   =\frac12  \int_{\mathbb{T}^d}  \big[|\u_0|^2+|\H_0|^2\big].
\end{equation}

For the system \eqref{ep8}, Proposition \ref{P1} implies that
$$
\int_{\mathbb{T}^d}\mathcal{Q}_1(\u, \mathbf{V}^0) \mathbf{V}^0=0,\ \
\int_{\mathbb{T}^d}\mathcal{Q}_2(\mathbf{V}^0, \mathbf{V}^0) \mathbf{V}^0=0,
$$
from which the following equality follows.
\begin{align}\label{cfh}
\frac12 \int_{\mathbb{T}^d}| \mathbf{V}^0|^2+\frac{\theta}{2}  \int_{\mathbb{T}^d}|\nabla \mathbf{V}^0|^2
 = \frac12\int_{\mathbb{T}^d}| \mathbf{V}^0(t=0)|^2.
\end{align}

Using $\mathcal{L}_1(\frac{t}{\epsilon})\mathbf{V}^0$ as a test function and noticing
$\rho^\epsilon=1+\epsilon \varphi^\epsilon$, we obtain the following
weak formulation of the continuity equation \eqref{a2i}
\begin{align}\label{ch}
  \int_{\mathbb{T}^d}
\mathcal{L}_1\Big(\frac{t}{\epsilon}\Big)\mathbf{V}^0\varphi^\epsilon(t)
&+\frac{1}{\epsilon}\int^t_0\int_{\mathbb{T}^d}
\Big[\text{div}\Big(\mathcal{L}_2\Big(\frac{\tau}{\epsilon}\Big)\mathbf{V}^0\Big)\varphi^\epsilon+\dv(\rho^\epsilon
\u^\epsilon)
  \mathcal{L}_1\Big(\frac{\tau}{\epsilon}\Big)\mathbf{V}^0\Big]\nonumber\\
 &   -\int^t_0\int_{\mathbb{T}^d}
 \mathcal{L}_1\Big(\frac{\tau}{\epsilon}\Big)\partial_t\mathbf{V}^0 \varphi^\epsilon
 =\int_{\mathbb{T}^d}
\varphi_0\varphi^\epsilon_0.
\end{align}
We use $\u$ and $\mathcal{L}_2 (\frac{t}{\epsilon} )\mathbf{V}^0$ to
test the momentum equation \eqref{a2j} respectively, to deduce
\begin{align}\label{civ}
&  \int_{\mathbb{T}^d}(\rho^\epsilon \u^\epsilon \cdot \u)(t)
+\int^t_0\int_{\mathbb{T}^d}\rho^\epsilon \u^\epsilon \cdot \big[(\u\cdot \nabla)
\u-(\H\cdot \nabla)\H-\mu \Delta \u+\nabla p+\frac12\nabla(|\H|^2)\big] \nonumber\\
&  -\int^t_0\int_{\mathbb{T}^d}\big[(\rho^\epsilon \u^\epsilon \otimes
\u^\epsilon)\cdot \nabla \u+   (\H^\epsilon\cdot \nabla)\H^\epsilon \cdot \u
 - \mu   \nabla  \u^\epsilon  \cdot \nabla \u\big] = \int_{\mathbb{T}^d}
  \rho^\epsilon_0 \u^\epsilon_0 \cdot \u_0
\end{align}
and
\begin{align}\label{cj}
&  \int_{\mathbb{T}^d}\Big(\rho^\epsilon \u^\epsilon \cdot \mathcal{L}_2
\Big(\frac{t}{\epsilon}\Big)\mathbf{V}^0\Big)(t)
+\frac{1}{\epsilon}\int^t_0\int_{\mathbb{T}^d}\rho^\epsilon \u^\epsilon
\cdot\nabla\Big(
\mathcal{L}_1\Big(\frac{\tau}{\epsilon}\Big)\mathbf{V}^0\Big)\nonumber\\
&\quad -\int^t_0\int_{\mathbb{T}^d}
\mathcal{L}_2\Big(\frac{\tau}{\epsilon}\Big)\partial_t\mathbf{V}^0
\cdot(\rho^\epsilon \u^\epsilon)
-\int^t_0\int_{\mathbb{T}^d}(\rho^\epsilon \u^\epsilon \otimes
\u^\epsilon)\cdot
\nabla\Big( \mathcal{L}_2\Big(\frac{\tau}{\epsilon}\Big)\mathbf{V}^0\Big)\nonumber\\
& \quad  +\int^t_0\int_{\mathbb{T}^d}\Big[\mu \nabla \u^\epsilon
\cdot \nabla \Big(
\mathcal{L}_2\Big(\frac{\tau}{\epsilon}\Big)\mathbf{V}^0\Big)
+(\mu+\lambda)\text{div}\u^\epsilon \,\text{div}
 \Big( \mathcal{L}_2\Big(\frac{\tau}{\epsilon}\Big)\mathbf{V}^0\Big)\Big]\nonumber\\
& \quad - \int^t_0\int_{\mathbb{T}^d}(\H^\epsilon\cdot
\nabla)\H^\epsilon \cdot
\mathcal{L}_2\Big(\frac{\tau}{\epsilon}\Big)\mathbf{V}^0
 -\int^t_0\int_{\mathbb{T}^d}\frac12|\H^\epsilon|^2\dv\Big(\mathcal{L}_2
 \Big(\frac{\tau}{\epsilon}\Big)\mathbf{V}^0\Big)\nonumber\\
& \quad
 -\int^t_0\int_{\mathbb{T}^d}
\Big(\frac{1}{\epsilon}\varphi^\epsilon+\frac{\gamma-1}{2}(\Pi^\epsilon)^2\Big)\text{div}\,
\Big(\mathcal{L}_2\Big(\frac{\tau}{\epsilon}\Big)\mathbf{V}^0\Big)
= \int_{\mathbb{T}^d}\rho^\epsilon_0 \u^\epsilon_0 \cdot
Q\tilde{\u}_0.
\end{align}
Similarly, we test \eqref{a2k} by $\H$ to get
\begin{align}\label{cjj}
&  \int_{\mathbb{T}^d}(\H^\epsilon \cdot \H)(t)
+\int^t_0\int_{\mathbb{T}^d}\H^\epsilon \cdot\big[(\u\cdot \nabla)
\H-(\H\cdot \nabla)\u-\nu \Delta \H\big] + \nu
\int^t_0\int_{\mathbb{T}^d}
\nabla \H^\epsilon  \cdot \nabla \H\nonumber\\
& \quad +\int^t_0\int_{\mathbb{T}^d}\big[(\dv  {\u}^\epsilon) { \H}^\epsilon +
( {\u}^\epsilon \cdot \nabla) {\H}^\epsilon  - ( {\H}^\epsilon\cdot \nabla)
 {\u}^\epsilon\big]\cdot \H  = \int_{\mathbb{T}^d}\H^\epsilon_0 \cdot \H_0.
\end{align}
Summing up \eqref{ce}, \eqref{cff} and \eqref{cfh},   inserting
\eqref{ch}-\eqref{cjj} into the resulting inequality,
and using the fact
\begin{align*}
 \int_0^t \int_{\mathbb{T}^d} \Big(\mathcal{L}\Big(\frac{\tau}{\epsilon}
 \Big)\partial_t\mathbf{V}^0\Big) \cdot \mathbf{U}^\epsilon
 =\int_0^t \int_{\mathbb{T}^d} \partial_t\mathbf{V}^0  \cdot \mathbf{V}^\epsilon,
\end{align*}
 we deduce, after a straightforward calculation, the following inequality:
\begin{align}\label{ckkl}
&   \frac{1}{2}\int_{\mathbb{T}^d}\Big\{  \Big|\sqrt{\rho^\epsilon} \u^\epsilon
-\u- \mathcal{L}_2\Big(\frac{t}{\epsilon}\Big)\mathbf{V}^0\Big|^2(t)+ |  \H^\epsilon
-\H|^2(t)
  +\Big|\Pi^\epsilon- \mathcal{L}_1\Big(\frac{t}{\epsilon}\Big)\mathbf{V}^0\Big|^2(t)\Big\}\nonumber\\
&\quad +  {\mu}  \int^t_0\int_{\mathbb{T}^d}\Big|\nabla\Big(
\u^\epsilon -\u-
\mathcal{L}_2\Big(\frac{\tau}{\epsilon}\Big)\mathbf{V}^0\Big)\Big|^2
+ \frac{\nu}{2}  \int^t_0\int_{\mathbb{T}^d}|\nabla \H^\epsilon-\nabla \H |^2 \nonumber\\
&\quad+ \frac{\mu}{2}  \int^t_0\int_{\mathbb{T}^d}\big(|\nabla \H^\epsilon|^2
+ |\nabla \H|^2\big)
+ (\mu+\lambda)\int^t_0\! \int_{\mathbb{T}^d}\Big|{\rm div}\Big( \u^\epsilon
-\u- \mathcal{L}_2\Big(\frac{\tau}{\epsilon}\Big)\mathbf{V}^0\Big)\Big|^2\nonumber\\
&\leq
  \frac{1}{2}\int_{\mathbb{T}^d}\Big\{  \Big|\sqrt{\rho^\epsilon} \u^\epsilon
-\u-
\mathcal{L}_2\Big(\frac{\tau}{\epsilon}\Big)\mathbf{V}^0\Big|^2(0)+
|  \H^\epsilon -\H|^2(0)
  +\Big|\Pi^\epsilon- \mathcal{L}_1\Big(\frac{\tau}{\epsilon}\Big)\mathbf{V}^0\Big|^2(0)\Big\}\nonumber\\
& \quad
+\sum^8_{i=1}A^\epsilon_i(t),
\end{align}
where
\begin{align}
A^\epsilon_1(t) =&
 \int_{\mathbb{T}^d}
\Big[(\sqrt{\rho^\epsilon}-1)\sqrt{\rho^\epsilon}\u^\epsilon\cdot
\Big(\u+\mathcal{L}_2\Big(\frac{t}{\epsilon}\Big)\mathbf{V}^0
\Big)\Big](t)-\int_{\mathbb{T}^d}
\Big[(\Pi^\epsilon-\varphi^\epsilon)\mathcal{L}_1\Big(\frac{t}{\epsilon}\Big)\mathbf{V}^0
\Big](t)\nonumber\\
& -\int_{\mathbb{T}^d}
\Big[(\sqrt{\rho^\epsilon}-1)\sqrt{\rho^\epsilon}\u^\epsilon\cdot
\Big(\u+\mathcal{L}_2\Big(\frac{t}{\epsilon}\Big)\mathbf{V}^0
\Big)\Big](0) +\int_{\mathbb{T}^d}
\Big[(\Pi^\epsilon-\varphi^\epsilon)\mathcal{L}_1\Big(\frac{t}{\epsilon}\Big)\mathbf{V}^0
\Big](0), \label{A1}\\
A^\epsilon_2(t) =&\int^t_0\int_{\mathbb{T}^d}\rho^\epsilon
\u^\epsilon \cdot
 \nabla p-\mu\int^t_0\int_{\mathbb{T}^d}(\rho^\epsilon-1)\u^\epsilon\Delta \u,\label{A2}\\
A^\epsilon_3 (t)= &
-\frac{\theta}{2}\int^t_0\int_{\mathbb{T}^d}\Delta  \mathbf{V}^0
\cdot \mathbf{V}^\epsilon
 - \mu\int^t_0\int_{\mathbb{T}^d}\nabla \u^\epsilon\cdot\nabla
 \Big( \mathcal{L}_2\Big(\frac{\tau}{\epsilon}\Big)\mathbf{V}^0\Big)\nonumber\\
&- (\mu+\lambda)\int^t_0\int_{\mathbb{T}^d}  \dv\u^\epsilon\cdot\dv
\Big( \mathcal{L}_2\Big(\frac{\tau}{\epsilon}\Big)\mathbf{V}^0\Big),\label{A3}\\
A^\epsilon_4 (t)= &-
\frac{\theta}{2}\int^t_0\int_{\mathbb{T}^d}|\nabla \mathbf{V}^0|^2
+\mu\int^t_0\int_{\mathbb{T}^d}\Big|\nabla
\mathcal{L}_2\Big(\frac{\tau}{\epsilon}\Big) \mathbf{V}^0\Big|^2
+(\mu+\lambda)\int^t_0\int_{\mathbb{T}^d}\Big|\dv \mathcal{L}_2
\Big(\frac{\tau}{\epsilon}\Big)\mathbf{V}^0\Big|^2,\label{A4}\\
A^\epsilon_5(t)=&  -\int^t_0\int_{\mathbb{T}^d}\rho^\epsilon
\u^\epsilon\cdot [(\H\cdot \nabla) \H)] -\int^t_0\int_{\mathbb{T}^d}
(\H^\epsilon\cdot \nabla) \H^\epsilon \cdot\u           \nonumber\\
&        +\int^t_0\int_{\mathbb{T}^d} \H^\epsilon \cdot\big[(\u\cdot \nabla)
\H-(\H\cdot \nabla)\H \big]  +\frac12\int^t_0\int_{\mathbb{T}^d}
\rho^\epsilon \u^\epsilon \cdot \nabla(|\H|^2)\nonumber\\
&   +\int^t_0\int_{\mathbb{T}^d}\big[(\dv  {\u}^\epsilon) { \H}^\epsilon
 + ( {\u}^\epsilon \cdot \nabla) {\H}^\epsilon
 - ( {\H}^\epsilon\cdot \nabla) {\u}^\epsilon\big]\cdot \H,  \label{A5}\\
A^\epsilon_6 (t)=& \int^t_0\int_{\mathbb{T}^d}\rho^\epsilon \u^\epsilon \cdot\big[ (\u\cdot
\nabla) \u   \big]-\int^t_0\int_{\mathbb{T}^d}(\rho^\epsilon
\u^\epsilon\otimes \u^\epsilon)\cdot \Big(\nabla u+  \mathcal{L}_2
\Big(\frac{\tau}{\epsilon}\Big)\mathbf{V}^0\Big) \nonumber\\
& - \frac{\gamma-1}{2}\int^t_0\int_{\mathbb{T}^d}
(\Pi^\epsilon)^2 \text{div} \Big( \mathcal{L}_2\Big(\frac{\tau}{\epsilon}\Big)\mathbf{V}^0\Big),\label{A6}\\
A^\epsilon_7 (t)= &  \int^t_0\int_{\mathbb{T}^d}\big[\mathcal{Q}_1(\u, \mathbf{V}^0)
  +\mathcal{Q}_2(\mathbf{V}^0, \mathbf{V}^0)\big]\cdot\mathbf{V}^\epsilon,\label{A7}\\
A^\epsilon_8 (t)= &
-\int^t_0\int_{\mathbb{T}^d}\frac12|\H^\epsilon|^2\dv\Big(\mathcal{L}_2
\Big(\frac{\tau}{\epsilon}\Big)\mathbf{V}^0\Big)
 - \int^t_0\int_{\mathbb{T}^d}(\H^\epsilon\cdot \nabla)\H^\epsilon
 \cdot \mathcal{L}_2\Big(\frac{\tau}{\epsilon}\Big)\mathbf{V}^0.\label{A8}
\end{align}

\medskip
\emph{Step 3: Convergence of the  modulated energy functional.}

To show the convergence of the modulated energy functional \eqref{ckkl},
 we need to estimate the remainders $A^{\epsilon}_i(t)$, $ i=1,\dots, 8$.
In the sequel, we will denote by $\omega^\epsilon(t)$ any sequence of time-dependent
functions which converges to $0$ uniformly in $t$.  For convenience, we also denote
 $\mathbf{w}^{\epsilon}\equiv\sqrt{\rho^\epsilon}
\u^\epsilon -\u-\mathcal{L}_2(\frac{t}{\epsilon})\mathbf{V}^0$,
$\mathbf{Z}^{\epsilon}\equiv \H^\epsilon -\H$, and $\Psi^\epsilon
\equiv \Pi^\epsilon-\mathcal{L}_1(\frac{t}{\epsilon})\mathbf{V}^0$.

For the  term $A^{\epsilon}_1(t)$, we employ \eqref{be},
\eqref{ine30},
 the regularity of  $\u$ and \eqref{UU},  and follow a procedure similar to
that for $R^\epsilon_2(t)$, to obtain
\begin{align}\label{dd1}
|A^{\epsilon}_1(t)|\leq C_T \epsilon +\omega^\epsilon(t).
\end{align}
On the other hand, the term $A^{\epsilon}_2(t)$ has the same bound as
$R^{\epsilon}_5(t)+R^{\epsilon}_7(t)$, thus
\begin{align}\label{dd2}
|A^{\epsilon}_2(t)|\leq C_T \epsilon.
\end{align}

To bound the term $A^{\epsilon}_3(t)$, we integrate by parts and use
the fact $ \mathcal{L}_2(\frac{t}{\epsilon})\mathbf{V}^0=\nabla
\tilde{q}^\epsilon$ for some function $\tilde{q}^\epsilon$  and
Proposition \ref{Pv} to infer
\begin{align*}
- \mu\int^t_0\int_{\mathbb{T}^d}\nabla \u^\epsilon\cdot\nabla\Big(
\mathcal{L}_2\Big(\frac{\tau}{\epsilon}\Big)\mathbf{V}^0\Big)
&=\mu\int^t_0\int_{\mathbb{T}^d}\Delta \u^\epsilon\cdot \mathcal{L}_2
\Big(\frac{\tau}{\epsilon}\Big)\mathbf{V}^0\\
&=\mu\int^t_0\int_{\mathbb{T}^d}\mathcal{L}\Big(-\frac{\tau}{\epsilon}\Big)
\Big(\begin{array}{c}
0\\
\Delta \u^\epsilon
\end{array}
\Big)\cdot  \mathbf{V}^0\\
&=\frac{\mu}{2}\int^t_0\int_{\mathbb{T}^d} \Delta \bar{\mathbf{V}}\cdot \mathbf{V}^0+\omega^\epsilon(t), \\
- (\mu+\lambda)\int^t_0\int_{\mathbb{T}^d}
\dv\u^\epsilon\cdot\dv\Big(
\mathcal{L}_2\Big(\frac{\tau}{\epsilon}\Big)\mathbf{V}^0\Big)
& =(\mu+\lambda)\int^t_0\int_{\mathbb{T}^d} \Delta \u^\epsilon\cdot\mathcal{L}_2
\Big(\frac{\tau}{\epsilon}\Big)\mathbf{V}^0\\
&=\frac{\mu+\lambda}{2}\int^t_0\int_{\mathbb{T}^d} \Delta \bar{\mathbf{V}}
\cdot \mathbf{V}^0+\omega^\epsilon(t),
\end{align*}
and
\begin{align*}
-\frac{\theta}{2}\int^t_0\int_{\mathbb{T}^d}\Delta  \mathbf{V}^0 \cdot \mathbf{V}^\epsilon
=-\frac{\theta}{2}\int^t_0\int_{\mathbb{T}^d}\Delta \bar{\mathbf{V}}\cdot \mathbf{V}^0+\omega^\epsilon(t).
\end{align*}
Thus, recalling $\theta=2\mu+\lambda$, one has
\begin{align}\label{dd3}
A^{\epsilon}_3(t)=\omega^\epsilon(t).
\end{align}

Similarly, it follows from Proposition \ref{Pv} that
\begin{align}\label{dd4}
A^{\epsilon}_4(t)=\omega^\epsilon(t).
\end{align}

Recalling \eqref{ck1i} and  using H\"older's inequality,  the
inequalities  \eqref{be} and \eqref{ine30},
    the regularity of  $\H$, and Sobolev's imbedding theorem, we conclude
\begin{align}\label{dd5}
A^{\epsilon}_5(t)
\leq & \int^t_0\int_{\mathbb{T}^d}(1-\rho^\epsilon)
\u^\epsilon\cdot [(\H\cdot \nabla) \H)]
+\int^t_0    ||\mathbf{Z}^\epsilon(\tau)||^2_{L^2}|| \nabla \u(\tau)||_{L^\infty}  d\tau\nonumber\\
&+  \int^t_0   \big[||\mathbf{w}^{\epsilon}(\tau)||^2_{L^2}
+ ||\mathbf{Z}^\epsilon(\tau)||^2_{L^2}\big] || \nabla \H(\tau)||_{L^\infty}  d\tau
\nonumber\\
&+\int^t_0\!\int_{\mathbb{T}^d}(\mathbf{Z}^\epsilon \cdot \nabla)\H \cdot
\Big[(1-\sqrt{\rho^\epsilon})\u^\epsilon+\mathcal{L}_2\Big(\frac{\tau}{\epsilon}\Big)\mathbf{V}^0\Big]
\nonumber\\
&- \int^t_0\!\int_{\mathbb{T}^d}\Big\{\Big[(1-\sqrt{\rho^\epsilon})\u^\epsilon
+\mathcal{L}_2\Big(\frac{\tau}{\epsilon}\Big)\mathbf{V}^0\Big]\cdot \nabla\Big\}
\H\cdot\mathbf{Z}^\epsilon\nonumber\\
& +\frac12\int^t_0\int_{\mathbb{T}^d}(\rho^\epsilon-1){\u}^\epsilon \nabla(|\H|^2)\nonumber\\
\leq &     \int^t_0   \big[||\mathbf{w}^{\epsilon}(\tau)||^2_{L^2}
+ ||\mathbf{Z}^\epsilon(\tau)||^2_{L^2}\big] \cdot \big[|| \nabla \u(\tau)||_{L^\infty}
+|| \nabla \H(\tau)||_{L^\infty}\big]  d\tau\nonumber\\
&+\int^t_0\!\int_{\mathbb{T}^d}\Big[(\mathbf{H}^\epsilon\cdot \nabla)\H \cdot
\mathcal{L}_2\Big(\frac{\tau}{\epsilon}\Big)\mathbf{V}^0 -
\Big(\mathcal{L}_2\Big(\frac{\tau}{\epsilon}\Big)\mathbf{V}^0\cdot
\nabla\Big)\H \cdot\mathbf{H}^\epsilon
 \Big]d\tau+C_T\epsilon.
\end{align}

The term $A^\epsilon_6(t)$ can be rewritten as
\begin{align}\label{dd6}
 A^\epsilon_6(t)   =& -\int^t_0\int_{\mathbb{T}^d}(\mathbf{w}^{\epsilon}\otimes
 \mathbf{w}^{\epsilon})\cdot \nabla\Big( \u +\mathcal{L}\Big(\frac{\tau}{\epsilon}\Big)\mathbf{V}^0\Big)
- \frac{\gamma-1}{2}\int^t_0\int_{\mathbb{T}^d}
|\Psi^\epsilon|^2 \text{div}\Big(\mathcal{L}_2\Big(\frac{\tau}{\epsilon}\Big)\mathbf{V}^0\Big) \nonumber\\
& -\int^t_0\int_{\mathbb{T}^d}\Big\{
\sqrt{\rho^\epsilon}\u^\epsilon\otimes \Big(\u
+\mathcal{L}_2\Big(\frac{\tau}{\epsilon}\Big)\mathbf{V}^0\Big)
+\Big(\u
+\mathcal{L}_2\Big(\frac{\tau}{\epsilon}\Big)\mathbf{V}^0\Big)\otimes
\sqrt{\rho^\epsilon}\u^\epsilon\Big\}
\cdot \nabla\Big( \u +\mathcal{L}\Big(\frac{\tau}{\epsilon}\Big)\mathbf{V}^0\Big)\nonumber\\
& +\int^t_0\int_{\mathbb{T}^d}
 \Big(\u +\mathcal{L}_2\Big(\frac{\tau}{\epsilon}\Big)\mathbf{V}^0\Big)
\otimes \Big(\u
+\mathcal{L}_2\Big(\frac{\tau}{\epsilon}\Big)\mathbf{V}^0\Big)\cdot
\nabla\Big( \u
+\mathcal{L}\Big(\frac{\tau}{\epsilon}\Big)\mathbf{V}^0\Big)
 \nonumber\\
     &+ \int^t_0\int_{\mathbb{T}^d}
\frac{\gamma-1}{2}
\Big|\mathcal{L}_1\Big(\frac{\tau}{\epsilon}\Big)\mathbf{V}^0\Big|^2
 \text{div} \Big( \mathcal{L}_2\Big(\frac{\tau}{\epsilon}\Big)\mathbf{V}^0\Big)
- {(\gamma-1)}
\mathcal{L}_1\Big(\frac{\tau}{\epsilon}\Big)\mathbf{V}^0\,
\Pi^\epsilon\,
\text{div} \Big( \mathcal{L}_2\Big(\frac{\tau}{\epsilon}\Big)\mathbf{V}^0\Big)\nonumber\\
&       +\int^t_0\int_{\mathbb{T}^d}(\rho^\epsilon-\sqrt{\rho^\epsilon})\u^\epsilon\cdot
 ((\u\cdot \nabla) \u)
-\int^t_0\int_{\mathbb{T}^d}(\sqrt{ \rho^\epsilon}\u^\epsilon- \u) \cdot \nabla
 \big(\frac{|\u|^2}{2}\big).
    \end{align}
We have to bound all the terms on the right-hand side of \eqref{dd6}.
Keeping in mind that $\dv\,\v=0$ and applying Proposition \ref{P2}, one obtains
\begin{align}\label{dd61}
& \int^t_0\int_{\mathbb{T}^d}
 \Big(\mathcal{L}_2\Big(\frac{\tau}{\epsilon}\Big)\mathbf{V}^0\Big)
\otimes
\Big(\mathcal{L}_2\Big(\frac{\tau}{\epsilon}\Big)\mathbf{V}^0\Big)\cdot
 \nabla\Big( \u +\mathcal{L}_2\Big(\frac{\tau}{\epsilon}\Big)\mathbf{V}^0\Big)
 +\frac{\gamma-1}{2}\Big|\mathcal{L}_1\Big(\frac{\tau}{\epsilon}\Big)\mathbf{V}^0\Big|^2
\dv\Big(\mathcal{L}_2\Big(\frac{\tau}{\epsilon}\Big)\mathbf{V}^0\Big)\nonumber\\
=& -\int^t_0\int_{\mathbb{T}^d}
 \Big\{\dv\Big(\mathcal{L}_2\Big(\frac{\tau}{\epsilon}\Big)\mathbf{V}^0\Big)
\otimes
\Big(\mathcal{L}_2\Big(\frac{\tau}{\epsilon}\Big)\mathbf{V}^0\Big)
+\frac{\gamma-1}{2}\nabla\Big|\mathcal{L}_1\Big(\frac{\tau}{\epsilon}\Big)
\mathbf{V}^0\Big|^2\Big\}
\cdot \big(\mathbf{V}^0 +\widehat{\u}\big)\nonumber\\
=& -\int^t_0\int_{\mathbb{T}^d} \mathcal{L}\Big(-\frac{\tau}{\epsilon}\Big)
\bigg(
\begin{array}{c}
0\\
\dv\Big(\mathcal{L}_2\Big(\frac{\tau}{\epsilon}\Big)\mathbf{V}^0\Big)
\otimes
\Big(\mathcal{L}_2\Big(\frac{\tau}{\epsilon}\Big)\mathbf{V}^0\Big)
 +\frac{\gamma-1}{2}\nabla\Big|\mathcal{L}_1\Big(\frac{\tau}{\epsilon}\Big)\mathbf{V}^0\Big|^2
\end{array}
\bigg)
\cdot \big(\mathbf{V}^0+\widehat{\u}\big)\nonumber\\
=&  -\int^t_0\int_{\mathbb{T}^d} \mathcal{Q}_2(\mathbf{V}^0,\mathbf{V}^0)
\cdot \big(\mathbf{V}^0+\widehat{\u}\big)+\omega^\epsilon(t)= \omega^\epsilon(t).
\end{align}
From Proposition \ref{P3} we get
\begin{align} \label{dd62}
&-\int^t_0\int_{\mathbb{T}^d}
\Big\{\mathcal{L}_2\Big(\frac{\tau}{\epsilon}\Big)\mathbf{V}^0\otimes\sqrt{\rho^\epsilon}\u^\epsilon
+
\sqrt{\rho^\epsilon}\u^\epsilon\otimes\mathcal{L}_2\Big(\frac{\tau}{\epsilon}\Big)\mathbf{V}^0\Big\}
\cdot\nabla\Big(\v +\mathcal{L}_2\Big(\frac{\tau}{\epsilon}\Big)\mathbf{V}^0\Big)\nonumber\\
& -\int^t_0\int_{\mathbb{T}^d}{(\gamma-1)}
\mathcal{L}_1\Big(\frac{\tau}{\epsilon}\Big)\mathbf{V}^0\,
\Pi^\epsilon\,
\text{div} \Big( \mathcal{L}_2\Big(\frac{\tau}{\epsilon}\Big)\mathbf{V}^0\Big)\nonumber\\
=&\int^t_0\int_{\mathbb{T}^d}
\Big\{\dv\Big(\mathcal{L}_2\Big(\frac{\tau}{\epsilon}\Big)\mathbf{V}^0\otimes\sqrt{\rho^\epsilon}\u^\epsilon
+\sqrt{\rho^\epsilon}\u^\epsilon\otimes\mathcal{L}_2\Big(\frac{\tau}{\epsilon}\Big)\mathbf{V}^0\Big)
+(\gamma-1) \nabla\Big(  \mathcal{L}_1\Big(\frac{\tau}{\epsilon}\Big)\mathbf{V}^0\Pi^\epsilon\Big)\Big\}\nonumber\\
& \ \ \ \ \quad  \cdot\Big(\u +\mathcal{L}_2\Big(\frac{\tau}{\epsilon}\Big)\mathbf{V}^0\Big)\nonumber\\
=&\int^t_0\int_{\mathbb{T}^d} \mathcal{L}\Big(-\frac{\tau}{\epsilon}\Big)
\bigg( \begin{array}{c}
0\\
\dv\Big(\mathcal{L}_2\Big(\frac{\tau}{\epsilon}\Big)\mathbf{V}^0\otimes\sqrt{\rho^\epsilon}\u^\epsilon
+
\sqrt{\rho^\epsilon}\u^\epsilon\otimes\mathcal{L}_2\Big(\frac{\tau}{\epsilon}\Big)\mathbf{V}^0\Big)
+(\gamma-1) \nabla\Big(
\mathcal{L}_1\Big(\frac{\tau}{\epsilon}\Big)\mathbf{V}^0\Pi^\epsilon\Big)
\end{array}
\bigg)\nonumber\\
& \ \ \ \ \quad \cdot \big(\mathbf{V}^0+\widehat{\u}\big)\nonumber\\
=&\int^t_0\int_{\mathbb{T}^d}
 \mathcal{L}\Big(-\frac{\tau}{\epsilon}\Big)
\bigg(
\begin{array}{c}
0\\
\dv\Big(\mathcal{L}_2\Big(\frac{\tau}{\epsilon}\Big)\mathbf{V}^0\otimes
Q(\sqrt{\rho^\epsilon}\u^\epsilon) +
Q(\sqrt{\rho^\epsilon}\u^\epsilon)\otimes\mathcal{L}_2\Big(\frac{\tau}{\epsilon}\Big)\mathbf{V}^0\Big)
+(\gamma-1) \nabla\Big(
\mathcal{L}_1\Big(\frac{\tau}{\epsilon}\Big)\mathbf{V}^0\Pi^\epsilon\Big)
\end{array}
\bigg)\nonumber\\
& \ \ \ \ \quad \cdot \big(\mathbf{V}^0+\widehat{\u}\big)\nonumber\\
&+\int^t_0\int_{\mathbb{T}^d}
 \mathcal{L}\Big(-\frac{\tau}{\epsilon}\Big)
\bigg(
\begin{array}{c}
0\\
\dv\Big(\mathcal{L}_2\Big(\frac{\tau}{\epsilon}\Big)\mathbf{V}^0\otimes
P(\sqrt{\rho^\epsilon}\u^\epsilon) +
P(\sqrt{\rho^\epsilon}\u^\epsilon)\otimes\mathcal{L}_2\Big(\frac{\tau}{\epsilon}\Big)\mathbf{V}^0\Big)
\end{array}
\bigg) \cdot \big(\mathbf{V}^0+\widehat{\u}\big)\nonumber\\
=& \int^t_0\int_{\mathbb{T}^d}\big\{ 2 \mathcal{Q}_2(\mathbf{V}^0, \bar{\mathbf{V}})
\cdot\big(\mathbf{V}^0+\widehat{\u}\big) +\mathcal{Q}_1(\u, \mathbf{V}^0)
\cdot\big(\mathbf{V}^0+\widehat{\u}\big)\big\}+ \omega^\epsilon(t).
\end{align}
Similarly, we have
\begin{align}\label{dd63}
&-\int^t_0\int_{\mathbb{T}^d}\big\{ \sqrt{\rho^\epsilon}\u^\epsilon\otimes \u
+\u\otimes \sqrt{\rho^\epsilon}\u^\epsilon\big\}
\cdot \nabla\Big( \u +\mathcal{L}\Big(\frac{\tau}{\epsilon}\Big)\mathbf{V}^0\Big)
=  \int^t_0\int_{\mathbb{T}^d}
 \mathcal{Q}_1(\u, \bar{\mathbf{V}}) \mathbf{V}^0
+ \omega^\epsilon(t),
\end{align}
and
\begin{align}\label{dd64}
& \int^t_0\int_{\mathbb{T}^d} \Big\{\u \otimes
\mathcal{L}_2\Big(\frac{\tau}{\epsilon}\Big)\mathbf{V}^0 +
\mathcal{L}_2\Big(\frac{\tau}{\epsilon}\Big)\mathbf{V}^0\Big)
\otimes \u \Big\} \cdot \nabla\Big( \u
+\mathcal{L}\Big(\frac{\tau}{\epsilon}\Big)\mathbf{V}^0\Big)
 \nonumber\\
& \quad = - \int^t_0\int_{\mathbb{T}^d}  \mathcal{Q}_1(\u, \mathbf{V}^0) \mathbf{V}^0
+ \omega^\epsilon(t).
\end{align}
On the other hand, using the basic energy equality \eqref{be} and the regularity of
$\u$, following the same process as for $R^\epsilon_6(t)$, we obtain
\begin{align}\label{dd65}
 \Big|       \int^t_0\int_{\mathbb{T}^d}(\rho^\epsilon-\sqrt{\rho^\epsilon})\u^\epsilon\cdot
 ((\u\cdot \nabla) \u)
-\int^t_0\int_{\mathbb{T}^d}(\sqrt{ \rho^\epsilon}\u^\epsilon- \u) \cdot \nabla
 \big(\frac{|\u|^2}{2}\big)\Big|
\leq C_T \epsilon.
\end{align}

The term $A^\epsilon_7(t)$ can be rewritten as
\begin{align}\label{dd7}
 A^\epsilon_7(t)
=  \int^t_0\int_{\mathbb{T}^d}\big[\mathcal{Q}_1(\u, \mathbf{V}^0)
  +\mathcal{Q}_2(\mathbf{V}^0, \mathbf{V}^0)\big]\cdot\bar{\mathbf{V}}+  \omega^\epsilon(t).
\end{align}
Substituting \eqref{dd61}-\eqref{dd7} into \eqref{dd6}, we conclude that
\begin{align}\label{A6A7}
 |A^\epsilon_6(t)| +|A^\epsilon_7(t)|\leq C\int^t_0(||\mathbf{w}^\epsilon||^2
 +||\Psi^\epsilon||^2)\Big(\|\nabla\u\|_{L^\infty} +\Big\|
\nabla\mathcal{L}_2\Big(\frac{\tau}{\epsilon}\Big)\mathbf{V}^0\Big\|_{L^\infty}\Big)d\tau+C_T
\epsilon + \omega^\epsilon(t).
\end{align}

Thus, we insert the estimates on $A^\epsilon_i(t)$ ($i=1,\cdots,7$)
into \eqref{ckkl} to obtain
\begin{align}\label{EE}
&   \frac{1}{2}\int_{\mathbb{T}^d}\big\{ | \mathbf{w}^\epsilon|^2
+ |   \mathbf{Z}^\epsilon  |^2
  + |\Psi^\epsilon|^2 \big\}(t)
+  {\mu}  \int^t_0\int_{\mathbb{T}^d}\Big|\nabla\Big( \u^\epsilon
-\u- \mathcal{L}_2\Big(\frac{\tau}{\epsilon}\Big)\mathbf{V}^0\Big)\Big|^2\nonumber\\
&\quad+ \frac{\nu}{2} \int^t_0\int_{\mathbb{T}^d} \big(|\nabla \mathbf{Z}^\epsilon |^2
+|\nabla \H^\epsilon|^2 + |\nabla \H|^2\big)
+ (\mu+\lambda)\int^t_0\! \int_{\mathbb{T}^d}\Big|{\rm div}\Big( \u^\epsilon
-\u- \mathcal{L}_2\Big(\frac{\tau}{\epsilon}\Big)\mathbf{V}^0\Big)\Big|^2\nonumber\\
\leq & \int^t_0 \big[||\mathbf{w}^{\epsilon}(\tau)||^2_{L^2}
+ ||\mathbf{Z}^\epsilon(\tau)||^2_{L^2}\big] \cdot \Big[|| \nabla
\u(\tau)||_{L^\infty}+|| \nabla \H(\tau)||_{L^\infty}
+\| \nabla\mathcal{L}_2\Big(\frac{\tau}{\epsilon}\Big)\mathbf{V}^0\|_{L^\infty})\Big]
d\tau  \nonumber\\
&  +\ \frac{1}{2}\int_{\mathbb{T}^d}\big\{ | \mathbf{w}^\epsilon|^2
+ |   \mathbf{Z}^\epsilon  |^2
  + |\Psi^\epsilon|^2\big\}(0)+C_T\epsilon +\omega^\epsilon(t)+
 A^\epsilon_8(t)+A^\epsilon_9(t),
\end{align}
where
\begin{align}\label{A9}
A^\epsilon_9(t)
=\int^t_0\!\int_{\mathbb{T}^d}\Big[(\mathbf{H}^\epsilon \cdot
\nabla)\H \cdot
\mathcal{L}_2\Big(\frac{\tau}{\epsilon}\Big)\mathbf{V}^0 -
\Big(\mathcal{L}_2\Big(\frac{\tau}{\epsilon}\Big)\mathbf{V}^0\cdot
\nabla\Big)\H \cdot\mathbf{H}^\epsilon
 \Big].
\end{align}

Now, to deal with $A^\epsilon_8(t)$ and $A^\epsilon_9(t)$, we denote
${\bar{\H}}^\epsilon_0=\frac{1}{|\mathbb{T}^d|}\int_{\mathbb{T}^d} \H^\epsilon_0(x)dx$ to
deduce from the magnetic field equation \eqref{a2k} that
\begin{equation}\label{average2}
\int_{\mathbb{T}^d} \H^\epsilon(x,t)dx= \int_{\mathbb{T}^d}
\H^\epsilon_0(x)dx = {\bar{\H}}^\epsilon_0 |\mathbb{T}^d|.
\end{equation}
The assumption that $\H^\epsilon_0$ converges strongly in $L^2$ to some
$\H_0$ implies
\begin{align*}
\Big|\int_{\mathbb{T}^d} \H^\epsilon_0(x)-\H_0(x)dx\Big|
 & \leq  \int_{\mathbb{T}^d} |\H^\epsilon_0(x)-\H_0(x)|dx   \\
& \leq |\mathbb{T}^d|^{\frac12} \Big( \int_{\mathbb{T}^d}|
\H^\epsilon_0(x)-\H_0(x)|^2dx\Big)^{\frac12} \rightarrow 0\ \
\text{as} \ \ \epsilon \to 0,
\end{align*}
whence
\begin{align}\label{average1}
{\bar{\H}}^\epsilon_0 \rightarrow {\bar{\H}}_0  \ \ \text{as} \ \
\epsilon\rightarrow 0.
\end{align}

Using H\"older's inequality, Sobolev's imbedding theorem, Poinc\'{a}re's
inequality, the isometry property of $\mathcal{L}$, and
\eqref{average2} and \eqref{average1}, we find that
\begin{align}\label{A89}
|{A}^\epsilon_8 (t)|+|{A}^\epsilon_9 (t)|
\leq &\frac{\nu}{4} \int^t_0[ ||\nabla\H^\epsilon||^2+ ||\nabla\H^\epsilon||^2](\tau)
d\tau+\omega(t)\nonumber\\
&+\frac{1}{\nu} (\|\varphi_0\|^2_{H^2}+||Q\u_0||^2_{H^2}) \int^t_0[
||\nabla\H^\epsilon||^2+ ||\nabla\H^\epsilon||^2](\tau)d\tau .
\end{align}

Thus, substituting \eqref{A89} into \eqref{EE} and using \eqref{Qu},
we deduce by Gronwall's inequality that, for almost all $t\in
[0,T]$,
\begin{align}\label{cqq}
 &
 ||\mathbf{w}^{\epsilon}(t)||^2_{L^2}+
 ||\mathbf{Z}^{\epsilon}(t)||^2_{L^2}+||\Psi^\epsilon(t)||^2_{L^2}\nonumber\\
& \quad \leq \tilde{C}
\Big\{||\mathbf{w}^{\epsilon}(0)||^2_{L^2}+||\mathbf{Z}^{\epsilon}(0)||^2_{L^2}
+ ||\Pi^\epsilon_0-\varphi_0||^2_{L^2} +C_T\epsilon
+\sup_{0\leq s\leq t}\omega^{\epsilon}(s)\Big\},
\end{align}
where
$$  \tilde{C}=\exp {\Big\{C\int^T_0 \Big[|| \nabla
\u(\tau)||_{L^\infty}+|| \nabla \H(\tau)||_{L^\infty}
+\| \nabla\mathcal{L}_2\Big(\frac{\tau}{\epsilon}\Big)\mathbf{V}^0\|_{L^\infty})\Big]
d\tau\Big\}}<+\infty. $$

\medskip
\emph{Step 4: End of the proof of Theorem \ref{MRc}.}

Letting $\epsilon$ go to zero in \eqref{cqq}, we see that
$\H^\epsilon$ converges strongly to
 $\H$ in $L^\infty(0,T; L^2(\mathbb{T}^d))$. Hence, $\bar{\mathbf{K}}=\mathbf{H}$.
Combining  \eqref{EE} with \eqref{cqq}, we can easily prove that
 $\nabla \H^\epsilon$ converges strongly to $\nabla \H$
in $L^2 (0,T;L^2(\mathbb{T}^d))$.

Next, it suffices to prove (4) and (5) in Theorem \ref{MRc}.
Noting that $P(\mathcal{L}_2(\frac{t}{\epsilon})\mathbf{V}^0)=0$ and the
fact that the projection operator $P$ is a bounded linear mapping form $L^2$ to $L^2$,
we obtain, with the help of \eqref{cqq}, that
\begin{align}\label{pes1}
\sup_{0\leq t\leq T}\|P(\sqrt{\rho^\epsilon} \u^\epsilon)
-\u\|_{L^2}
=&\sup_{0\leq t\leq T}\Big\|P\Big(\sqrt{\rho^\epsilon} \u^\epsilon-\u-\mathcal{L}_2
\Big(\frac{t}{\epsilon}\Big)\mathbf{V}^0\Big)\Big\|_{L^2}\nonumber\\
\leq &\sup_{0\leq t\leq T}\Big\|\sqrt{\rho^\epsilon} \u^\epsilon-\u-\mathcal{L}_2
\Big(\frac{t}{\epsilon}\Big)\mathbf{V}^0\Big\|_{L^2}\nonumber\\
\rightarrow&\;\; 0 \;\;\mbox{as}\;\; \epsilon \rightarrow 0.
\end{align}
Therefore, (4) is proved. Utilizing \eqref{rhoc}, we deduce from \eqref{a2i}
that $\mbox{div}\,\u^\epsilon$ converges weakly to $0$ in $H^{-1}((0, T)\times\mathbb{T}^d)$.
Thus we obtain easily that $Q\u^\epsilon$ converges weakly to $0$
in $H^{-1}(0, T; L^2(\mathbb{T}^d))$. In view of the fact that
$\u^\epsilon$ is bounded in $L^2(0,T;L^2(\mathbb{T}^d))$ and $\sqrt{\rho^\epsilon}$
converges strongly to $1$ in $C([0, T], L^2(\mathbb{T}^d))$, we see that
$Q(\sqrt{\rho^\epsilon}\u^\epsilon)$ converges weakly to $0$ in
$H^{-1}(0, T; L^2(\mathbb{T}^d))$. Obviously, the fact $\sqrt{\rho^\epsilon}\u^\epsilon=
P(\sqrt{\rho^\epsilon} \u^\epsilon)+Q(\sqrt{\rho^\epsilon}\u^\epsilon)$
implies the weak convergence of $\sqrt{\rho^\epsilon}\u^\epsilon$
to $\u$ in $H^{-1}(0, T; L^2(\mathbb{T}^d))$. The proof of Theorem \ref{MRc} is finished.
\end{proof}
Theorem \ref{MRd} can be shown by slightly modifying the proof of Theorem \ref{MRc}, and
therefore, we omit its proof here.
\medskip

{\bf Acknowledgement:} This work was done while Fucai Li was
visiting the Institute of Applied Physics and Computational
Mathematics. He would like to thank the institute for hospitality.
Jiang was supported by the National Basic Research Program (Grant
No. 2005CB321700) and NSFC (Grant No. 40890154). Ju was supported by
NSFC (Grant No. 10701011). Li was supported  by NSFC (Grant No.
10501047, 10971094).


\end{document}